\documentclass[11pt]{article}

\usepackage{graphicx}

\usepackage{amsthm}
\usepackage{amsmath}
\usepackage{amsfonts}
\usepackage{amssymb}
\usepackage{array}
\usepackage{bm}
\usepackage{multirow}
\usepackage[shortlabels, inline]{enumitem}
\usepackage[normalem]{ulem}
\usepackage[font=small]{caption}
\usepackage[a4paper]{geometry}
\usepackage{algorithm}
\usepackage{algpseudocode}

\algnewcommand{\IIf}[1]{\State\algorithmicif\ #1\ \algorithmicthen}

\usepackage[round]{natbib}
\usepackage[colorlinks=true, linkcolor=red, citecolor=blue]{hyperref}

\newtheorem{theorem}{Theorem}[section]
\newtheorem{lemma}[theorem]{Lemma}
\newtheorem{corollary}[theorem]{Corollary}
\newtheorem{proposition}[theorem]{Proposition}

\newtheorem*{theorem*}{Theorem}
\newtheorem*{lemma*}{Lemma}
\newtheorem*{corollary*}{Corollary}
\newtheorem*{proposition*}{Proposition}
\newtheorem*{conjecture*}{Conjecture}

\theoremstyle{definition}
\newtheorem{definition}{Definition}[section]
\newtheorem*{definition*}{Definition}
\theoremstyle{definition}

\theoremstyle{definition}
\newtheorem{example}{Example}
\theoremstyle{remark}
\newtheorem*{example*}{Example}
\theoremstyle{definition}

\theoremstyle{definition}
\newtheorem*{assm*}{Assumption}
\theoremstyle{remark}
\newtheorem{remark}{Remark}[section]
\theoremstyle{remark}
\newtheorem*{remark*}{Remark}

\def\filedate{4/26/2017}
\def\fileversion{v0.2}
\RequirePackage{kvoptions}
\SetupKeyvalOptions{
    family=@arash,
    prefix=@arash@
}
\DeclareBoolOption{matrixnorms}\DeclareStringOption{opnorm}
\ProcessKeyvalOptions*

\newcommand{\reals}{\mathbb{R}}

\newcommand{\nats}{\mathbb{N}}

\newcommand{\pr}{\mathbb{P}}
\newcommand{\eps}{\varepsilon}
\newcommand{\ex}{\mathbb{E}}

\DeclareMathOperator{\supp}{supp}

\newcommand\vnorm[1]{\|#1\|}
\newcommand\norm[1]{\|#1\|}

\if@arash@matrixnorms
    \newcommand\mnorm[1]{|\!|\!|#1|\!|\!|}
\else
    \newcommand\mnorm[1]{\vnorm{#1}}    
\fi

\newcommand\opnorm[1]{\mnorm{#1}_{\text{\@arash@opnorm}}}

\newcommand{\ran}{ran}

\RequirePackage{tikz}

\newcommand{\ip}[1]{\langle #1\rangle}

\makeatletter
\let\save@mathaccent\mathaccent
\newcommand*\if@single[3]{\setbox0\hbox{${\mathaccent"0362{#1}}^H$}\setbox2\hbox{${\mathaccent"0362{\kern0pt#1}}^H$}\ifdim\ht0=\ht2 #3\else #2\fi
}
\newcommand*\rel@kern[1]{\kern#1\dimexpr\macc@kerna}
\newcommand*\widebar[1]{\@ifnextchar^{{\wide@bar{#1}{0}}}{\wide@bar{#1}{1}}}
\newcommand*\wide@bar[2]{\if@single{#1}{\wide@bar@{#1}{#2}{1}}{\wide@bar@{#1}{#2}{2}}}
\newcommand*\wide@bar@[3]{\begingroup
    \def\mathaccent##1##2{\let\mathaccent\save@mathaccent
\if#32 \let\macc@nucleus\first@char \fi
\setbox\z@\hbox{$\macc@style{\macc@nucleus}_{}$}\setbox\tw@\hbox{$\macc@style{\macc@nucleus}{}_{}$}\dimen@\wd\tw@
        \advance\dimen@-\wd\z@
\divide\dimen@ 3
        \@tempdima\wd\tw@
        \advance\@tempdima-\scriptspace
\divide\@tempdima 10
        \advance\dimen@-\@tempdima
\ifdim\dimen@>\z@ \dimen@0pt\fi
\rel@kern{0.6}\kern-\dimen@
        \if#31
        \overline{\rel@kern{-0.6}\kern\dimen@\macc@nucleus\rel@kern{0.4}\kern\dimen@}\advance\dimen@0.4\dimexpr\macc@kerna
\let\final@kern#2\ifdim\dimen@<\z@ \let\final@kern1\fi
        \if\final@kern1 \kern-\dimen@\fi
        \else
        \overline{\rel@kern{-0.6}\kern\dimen@#1}\fi
    }\macc@depth\@ne
    \let\math@bgroup\@empty \let\math@egroup\macc@set@skewchar
    \mathsurround\z@ \frozen@everymath{\mathgroup\macc@group\relax}\macc@set@skewchar\relax
    \let\mathaccentV\macc@nested@a
\if#31
    \macc@nested@a\relax111{#1}\else
\def\gobble@till@marker##1\endmarker{}\futurelet\first@char\gobble@till@marker#1\endmarker
    \ifcat\noexpand\first@char A\else
    \def\first@char{}\fi
    \macc@nested@a\relax111{\first@char}\fi
    \endgroup
}
\makeatother

\usepackage{pgf, tikz}
\usetikzlibrary{arrows, automata}
\usepackage{subcaption}

\providecommand{\Xc}{\mathcal{X}}
\providecommand{\Tc}{\mathcal{T}}

\providecommand{\Hil}{\mathcal{H}}
\DeclareMathOperator{\Span}{span}
\def\SEP(#1,#2,#3){#1-#2-#3}

\providecommand{\Ss}{S^*}

\providecommand{\nhbdcoef}{\beta}

\providecommand{\nhbdset}{\mathcal{T}}

\providecommand{\dvec}{X}

\providecommand{\posdef}{\succ}

\providecommand{\normalN}{\mathcal{N}}

\providecommand{\coef}{\beta}
\providecommand{\R}{\mathbb{R}} \providecommand{\reals}{\mathbb{R}}
 \providecommand{\ex}{\mathbb{E}}

\providecommand{\tpose}{T}

\providecommand\independent{\protect\mathpalette{\protect\independenT}{\perp}}
\def\independenT#1#2{\mathrel{\rlap{$#1#2$}\mkern4mu{#1#2}}}

\def\given{\,\mid\,}

\providecommand{\norm}[1]{\Vert#1\Vert}

\definecolor{maroon}{rgb}{0.5, 0.0, 0.0}

\providecommand{\Cc}{\mathcal{C}}
\newcommand{\Pc}{\mathcal{P}}

\providecommand{\rv}{X}
\providecommand{\gv}{X}

\newcommand{\Fc}{\mathcal F}

\usepackage{centernot}

\newcommand{\note}[1]{\textcolor{blue}{\textbf{\texttt{[#1]}}}}
\newcommand{\update}[1]{\textcolor{blue}{#1}}

\newcommand{\rev}[1]{\textcolor{teal}{#1}}
\newcommand{\revv}[1]{\textcolor{teal}{#1}}
\newcommand{\aarev}[1]{\textcolor{black!40!blue}{#1}}
\renewcommand{\rev}[1]{\textcolor{black}{#1}}
\renewcommand{\revv}[1]{\textcolor{black}{#1}}

\renewcommand{\note}[1]{}
\renewcommand{\update}[1]{\textcolor{black}{#1}}
\renewcommand{\aarev}[1]{\textcolor{black}{#1}}

\DeclareMathOperator{\cov}{cov}

\newcommand\res{\delta}
\newcommand\rhoh{\widehat\rho}

\newcommand\Fcb{\bar \Fc}
\newcommand\gcm{\varphi}
\newcommand\gcmh{\widehat \gcm}

\newcommand\gcmodel{\operatorname{GPC}}

\newcommand\gcorr{\rho}
\newcommand\gcorrh{\widehat \gcorr}

\newcommand\fh{\hat f}

\newcommand{\AND}{\text{ and }}

\newcommand{\nhbdlat}{\mathfrak{T}}
\newcommand{\nhbdlatred}{\mathfrak{t}}
\newcommand{\nhbdlatd}{\mathfrak{D}}

\newcommand{\ground}{V}

\newcommand{\prob}{\nu}
\newcommand{\gr}{G}
 \newcommand{\latmin}{m}
\newcommand{\latmax}{M}
\newcommand{\ind}{\mathcal{I}}
\newcommand{\graphoid}{\ind}

\newcommand{\indep}{\independent}
\newcommand{\notindep}{\centernot\indep}
\newcommand{\indepprob}{\indep_{\!\!\!\prob}\,}
\newcommand{\indepgr}{\indep_{\!\!\!\gr}\,}

\newcommand{\Bf}{\mathfrak{B}}
\newcommand\Bfin{\Bf_{\text{fin}}}
\DeclareMathOperator{\csp}{\overline{Sp}}

\DeclareMathOperator\Prj{Proj}
\newcommand\Hilset{\Gamma}

\usepackage{chngcntr}
\newcounter{exletter}
\newcounter{exnumber}
\counterwithin{exnumber}{exletter}


\title{Learning general conditional independence structures\\via the neighbourhood lattice
}
\author{Arash A. Amini, Bryon Aragam, Qing Zhou}

\begin{document}




\maketitle

\begin{abstract}We study the problem of learning multivariate dependencies in nonparametric and high-dimensional settings.
\update{This includes but is not limited to graphical models.}
Our approach effectively combines several features that are missing from previous work on this problem: We show how the entire dependence structure can be learned nonparametrically while simultaneously evading the curse of dimensionality and relaxing common assumptions such as faithfulness.
To this end, we introduce and study the neighbourhood lattice decomposition of a distribution, which is a compact, non-graphical representation of conditional independence (CI) that is valid in the absence of a faithful graphical representation.
We show that the neighbourhood lattice decomposition exists in any graphical model and can be computed efficiently, nonparametrically, and consistently in high-dimensions without paying the usual curse of dimensionality.
This gives a way to \rev{learn} \emph{all} of the independence relations implied by \emph{any} graphical model, without requiring \emph{a priori} knowledge of the graph or even the graph type.
As a special case, our results provide a general solution to the problem of nonparametric estimation of high-dimensional CI structures over any graphical model.
\end{abstract}


\section{Introduction}
\label{sec:intro}

We consider the problem of learning general dependence structures in high-dimensional and nonparametric models.
Classical special cases of this problem \rev{include support recovery in regression (i.e. learning a Markov boundary), structure learning in graphical models, and estimating causal effects from observational studies.}
Knowledge of the (in)dependence structure of a system can be leveraged
to perform efficient inference,
to represent and manipulate distributions in memory,
to express causal relationships,
to design experiments and \rev{select} interventions,
among many other fundamental tasks \citep{wainwright2008graphical,pearl2009,spirtes2000,koller2009,lauritzen1996}.

\rev{Representing independence structures via graphical models has been exploited with great success: Graphical models such as Markov random fields (MRFs) and Bayesian networks (BNs) allow for tractable inference and structure learning in special cases.
However, these models face inherent limitations.} It is well-known that even arbitrarily complex graphical models cannot represent all possible independence structures \citep{studeny2006probabilistic,sadeghi2017faithfulness}.
As a consequence, graphical models suffer from an intrinsic model misspecification issue: We must know in advance the \emph{type} of graphical model to use (e.g., BNs, MRFs, or more exotic variants), and even then, the chosen type might not be capable of faithfully representing the true underlying dependencies.

Motivated by these limitations, this paper introduces the \emph{neighbourhood lattice}, a novel approach to learn dependence structures that simultaneously addresses several challenges. Our framework allows us to learn the CI structure of \emph{any} graphical model (e.g., MRFs, BNs, chain graphs, etc.) without suffering the curse of dimensionality and without prior knowledge of its specific type or assuming faithfulness. More powerfully, it extends to CI structures that \emph{cannot} be represented by any graph at all. The key to this generality is a non-graphical representation of CI that, as we will demonstrate, remains computationally and statistically tractable even in high-dimensional, nonparametric settings. To build intuition for this general concept, let us first explore its core idea in a familiar context.

\subsection{The Neighbourhood Lattice: An Intuitive Introduction via Regression}
\label{sec:intro:intuitive_lattice}

To grasp the core idea behind our approach, consider a familiar setting: Understanding which covariates explain a response variable in linear regression.
Let $X=(X_{1},\ldots,X_{d})$ be a random vector. \aarev{ Our goal is to understand the dependencies between the components of $X$. To this end, identify $X_j$ as $Y$ (the response) and $X_{{-j}}$ as $Z = (Z_1,\dots,Z_{d-1})$ (the covariates).}
If we regress $Y$ on all covariates $Z$, we obtain a set of coefficients, some of which may be zero. The non-zero coefficients identify an \emph{active set} $m_1 \subseteq Z$, meaning $Y \indep Z_i \given Z_{m_1}$ for any $i \notin m_1$.\footnote{This is always true when $X$ is Gaussian with nonsingular covariance, however, our main results apply to general, non-Gaussian distributions. Extending this beyond regression and Gaussianity is a key motivation, discussed in the next subsection.} This $m_1$ is the \emph{Markov boundary} of $Y$ relative to $Z$. 
\update{The remaining Markov blankets of $Y$ satisfy $Y \indep Z_i \given Z_{U}$ for any $i \notin U$, and there is always a largest Markov blanket $M_1$, i.e. the largest superset of $m_1$ within $Z$ such that $Y \indep Z_{i} \given Z_{M_1}$ for $i \notin U$.}

A typical regression analysis assesses this local CI relation. However, what if we vary the set of covariates?
For example, what if we exclude a variable from $m_1$ or even the entire set $m_1$ and try to explain $Y$ with the remaining variables? Or, is there another subset of $Z$, perhaps of similar size to $m_1$, that also renders the rest of $Z$ irrelevant for explaining $Y$?
Each time we do this, we learn more about the CI structure, however, answering these questions typically requires re-running the regression for each new candidate subset of covariates $S_\ell \subseteq Z$, yielding new active sets $m_\ell$ and maximal sets $M_\ell$.

Our key insight is that one does not need to perform all these regressions. There is significant overlap.
Using set interval notation, $$[m_\ell, M_\ell] := \{U \subseteq X_{-j}:\; m_\ell \subseteq U \subseteq M_{\ell} \},$$ any subset $U$ in this interval will have the same explanatory power for $Y$ as $m_\ell$.
Remarkably, the collection of distinct intervals $[m_\ell, M_\ell]$ forms a partition of the power set $2^{X_{-j}}$.
This partition, which we call the \emph{neighbourhood lattice decomposition} for node $j$, directly represents all local CI relations of the form $Y \indep B \given C$ for all choices of $B,C \subseteq Z$.
As we will show, this simple idea, when generalized beyond regression and Gaussianity, provides a powerful, non-graphical way to represent and learn \emph{all} conditional independencies in a distribution, even those that cannot be captured by any graphical model.

\subsection{A Key Implication: Circumventing the Curse of Dimensionality}
\label{sec:intro:payoff}

One powerful consequence of this neighbourhood lattice approach is the ability to learn these complex CI structures nonparametrically in high dimensions \emph{without succumbing to the curse of dimensionality}.
To illustrate, consider the challenge of evaluating a CI relation like $X_1 \indep X_2 \given X_{3:d}$ in a general $d$-dimensional distribution. This is comparable to estimating the regression function $\ex[X_1 \given X_{2:d}]$, a task notoriously hampered by the curse of dimensionality: For an $s$-smooth regression function, the minimax $L_2$ estimation rate is $n^{-s/(2s+d)}$, which degrades rapidly with $d$.

The neighbourhood lattice offers a path to overcome this. 
\update{Instead of suffering the curse of dimensionality, the effective dimension for nonparametric estimation of the CI structure is a sparsity parameter $t$ instead of the ambient dimension $d$.}
We will establish results of the following informal nature:
\begin{theorem}[Informal]\label{thm:informal_intro}
    Let $t \in [d]$ and $\alpha, s > 0$.
    Assume that $(\rv_{1},\ldots,\rv_{d})$ has a distribution $\nu>0$ supported on a bounded set, with a CI structure that is (a) a compositional graphoid (a general property satisfied by any graphical model, see Section~\ref{sec:bg}, Lemma~\ref{lem:ex:cg}), (b) $\alpha$-separated of order $t$ (a signal strength condition), and (c) $t$-sparse (meaning its Markov boundaries are of size at most $t$). Moreover, for any $|S| \le t$ assume that the regression functions $\ex[X_i \given X_S]$ are $s$-H\"{o}lder smooth.
    Then, assuming $(2-\eps)s \ge t$ for some positive $\eps$, if
    \begin{align}\label{eq:alpha:lower:bound:informal_intro}
    \alpha^2 \;\gtrsim\; n^{-2s / (2s+ t)} + \frac{t \log d}{n},
    \end{align}
    then all the CI relationships in $\nu$ can be correctly identified from a sample of size $n$, with high probability.
\end{theorem}
\noindent
This result (formalized in Theorem~\ref{thm:highd:consist:nongauss}) demonstrates that the effective dimension for nonparametric estimation becomes a sparsity parameter $t$, not the ambient dimension $d$. This allows consistent CI discovery even when $d \gg n$, including for CI relations with conditioning sets much larger than $t$, paying only a logarithmic factor in $d$ and the curse of dimensionality in $t$.
Thus, the neighbourhood lattice provides a general solution to nonparametric estimation of high-dimensional CI structures, opening fruitful avenues for future investigation.

\subsection{Our Contributions}
\label{sec:intro:contributions}

To achieve these goals, this paper develops the following:
\begin{enumerate}[itemsep=4pt]
\item \emph{Representing CI structures via the neighbourhood lattice} (Section~\ref{sec:lat}).
We detail the construction and elementary properties of the neighbourhood lattice (Theorem~\ref{thm:nhbdlat}) and the resulting neighbourhood lattice decomposition (Definition~\ref{defn:lat:decomp}). This decomposition represents the entire CI structure of any distribution $\prob$ that satisfies certain general axioms (specifically, the \emph{composition} and \emph{intersection} axioms of graphoids, which encompasses graphical models like MRFs and BNs, as well as non-graphical structures, often without requiring faithfulness assumptions; see Section~\ref{sec:bg}). Any CI relation implied by $\prob$ can be easily checked using the neighbourhood lattice (Corollaries~\ref{cor:elem:ci}-\ref{cor:gen:ci}).

\item \emph{Estimation and computation} (Section~\ref{sec:comp}). We propose practical algorithms to learn the neighbourhood lattice (Algorithms~\ref{alg:lattice:comp}-\ref{alg:lattice:decomp:sparse})\rev{---and hence the entire CI structure---}and analyze their computational complexity (Theorem~\ref{thm:poly:N}).
\rev{Open-source code implementing the main algorithms is available at~\url{https://github.com/aaamini/nblat/}.}

\item \emph{Consistency and finite-sample properties} (Section~\ref{sec:highdim}). We prove the consistency of our algorithms and study the sample complexity of estimating the neighbourhood lattice decomposition for general classes of distributions in high-dimensions (Theorem~\ref{thm:highd:consist:nongauss}). \rev{Our results indicate that the curse of dimensionality can be circumvented even when we seek to learn \emph{every} CI relation, including those with large conditioning sets.}

\item \rev{\emph{Connection to graphical models and regression} (Sections~\ref{sec:grint}-\ref{sec:proj}). We show by example how the neighbourhood lattice can be interpreted through well-known properties of graphical models, and also bears a close relationship with commonly studied regression models in statistics, especially the concept of \emph{neighbourhood regression} in graphical models. This illustrates how our framework subsumes graphical models as a special case.}
\end{enumerate}
For clarity and generality, although our main interest will be applications to probabilistic conditional independence, we will present our results in the setting of an abstract compositional graphoid (see Section~\ref{sec:bg} for definitions).
The rich algebraic structure of the neighbourhood lattice allows derivation of new CI relations ``for free" from a subset of known relations, significantly reducing the statistical and computational burden compared to brute-force approaches.

\subsection{Related Work}
\label{sec:intro:related}

The problem of representing \rev{and learning} probabilistic CI structures has a long history; we refer the reader to textbooks such as \cite{pearl1988,lauritzen1996,studeny2006probabilistic} for additional background.
Among the several representations of CI structures are
graphical models \citep{geiger1993,lauritzen1988local},
graphoids \citep{pearl1985graphoids,pearl1987logic},
separoids \citep{dawid2001separoids},
and imsets \citep{studeny1995description}.
\rev{Chain graphs \citep{lauritzen1989graphical, levitz2001separation}, maximal ancestral graphs \citep{richardson2002}, and nested Markov models \citep{richardson2023nested} represent efforts to generalize graphical models for richer CI structures and causal applications. \citet{lauritzen2018unifying} provided a conceptual unification, proving that all such models are compositional graphoids, as studied here. \citet{sadeghi2017faithfulness} established conditions for faithfulness to specific graphs.}
Generalizing beyond graphical models, the literature has introduced matroids, separoids, and imsets. 
However, learning and inference with these models is difficult, and there is no finite complete characterization of conditional independence \citep{studeny1990conditional}, even for restrictive families such as Gaussians \citep{sullivant2009gaussian}. Our neighbourhood lattice approach targets a middle ground, offering greater generality than graphical models while retaining computational and statistical tractability.

Lattice theory has been applied to CI relations before \citep[e.g.,][]{dawid2001separoids,niepert2013conditional,gaag2018lattice,andersson1993lattice, studeny2006probabilistic}, but the specific neighbourhood lattice we introduce and its decomposition appears to be novel.
Our work builds upon known facts about CI and structure learning, including Markov boundary learning. Constraint-based methods like the PC algorithm \citep{spirtes1991,kalisch07} are related, but rely on faithfulness assumptions. The GES algorithm works under compositionality, returning a minimal I-map \citep{chickering2002}, but this may not represent all CI relations. See also \citet{meek1997thesis,chickering2002ges}.
Recent work on the hardness of CI testing \citep{azadkia2019simple,canonne2018testing,neykov2020minimax,shah2020hardness} highlights the challenges, particularly in multivariate settings, which our work addresses by leveraging the lattice structure.
Finally, nonparametric approaches for specific graphical models like the nonparanormal \citep{liu2009} or subgaussian inverse covariance models \citep{ravikumar2011} exist, but our framework is substantially more general, subsuming these as special cases by applying to all compositional graphoids under only smoothness assumptions.

\subsection{Notation}
\label{sec:intro:notation}
For any index set $\ground$, $A\subseteq \ground$, and a random vector $(X_{j})_{j\in\ground}$, we denote $\rv_{A} = \{\rv_i:\; i \in A\}$.
For any integer $m$ let $[m]:=\{1,\ldots,m\}$; we will often take $d:=|\ground|$ and $V=[d]$.
We also use the shorthand notations: $\{i\} = i$ and $\{i,j\} = ij$, $A \cup \{i\} = Ai$, $A \cup B = AB$ and so on.
In addition, we let $\ground_{-S} = \ground \setminus S$. Common uses of these notational conventions are: $\ground_{-j} = \ground \setminus \{j\}$ and $\ground_{-ij} = \ground \setminus ij = \ground \setminus \{i,j\}$.
For a matrix $\Sigma \in \reals^{d \times d}$ and subsets $A,B \subseteq \ground$, we use $\Sigma_{A,B}$ for the submatrix on rows and columns indexed by $A$ and $B$, respectively. Single index notation is used for principal submatrices, so that $\Sigma_{A} = \Sigma_{A,\,A}$.
Finally, we use $A\uplus B$ to indicate the union of two disjoint sets $A,B\subseteq\ground$.
We write $X_{A}\indepprob X_{B}\given X_{C}$ for CI under distribution $\prob$, or $A\indep B \given C$ for an abstract independence model.

\section{Background}
\label{sec:bg}

To formally develop the neighbourhood lattice and prove its properties for a wide range of dependence structures—including all standard graphical models and beyond, often without assumptions like faithfulness—we employ the framework of abstract independence models \citep{lauritzen2018unifying,sadeghi2017faithfulness,studeny2006probabilistic}.
Although our results can be specialized to arbitrary graphical models, this more abstract setting is necessary to apply our results to non-graphical models and situations where faithfulness does not hold.
We will use the established framework of abstract (semi-)graphoids, of which probabilistic conditional independence and graphical models are special cases. This section collects the necessary background and definitions as well as several examples for motivation.
Readers less familiar with this abstraction can, on a first reading, ignore this abstraction and instead substitute their preferred graphical model (e.g., MRFs, BNs) and its associated CI properties without any modification to the main thrust of our results. Lemma~\ref{lem:ex:cg} later in this section confirms these are indeed special cases.

First, some preliminaries: Recall that a complete lattice is a partially ordered set (poset) in which all subsets have both a supremum (join) and an infimum (meet)~\cite[Section 3.3]{stanley1997enumerative}. We also need the following definition: Let $(\Pc, \le)$ be a poset and $\Cc$ a subposet of $\Pc$. We say that $\Cc$ is \emph{convex} (in $\Pc$) if $z \in \Cc$ whenever $x < z < y$ with $x,y \in \Cc$. An interval $[x,y] := \{z \in \Pc:\; x \le z \le y\}$ is an example of a convex subposet~\cite[Section~3.1]{stanley1997enumerative}.
Given a set $\Omega$, we are mostly interested in the poset $\Bf(\Omega)$, the set of all subsets of $\Omega$, ordered by inclusion, i.e., $\Bf(\Omega) = (2^\Omega, \subseteq)$. An interval $[A,B]$ in $\Bf(\Omega)$ consists of all subsets $S \subseteq \Omega$ such that $A \subseteq S \subseteq B$. It is easy to verify that $\Bf(\Omega)$ is a complete lattice with the join and meet given by the set union and set intersection, respectively, and any interval in $\Bf(\Omega)$ is a convex sublattice of $\Bf(\Omega)$.

\subsection{Graphoids and independence models}
\label{sec:gr}
Let $\rv$ be a random vector; we make no distributional or parametric assumptions on this vector such as Gaussianity, linearity, etc.
Let $V=[d]$ be the index set, so $\rv=(\rv_{j})_{j\in\ground}$.
As with graphical models, we abuse notation by identifying $\ground=\rv=[d]$.
A (formal) \emph{independence model} over $V$ is a ternary relation $A\indep B\given C$ over disjoint subsets $A,B,C\subseteq V$. A \emph{graphoid} over $V$ is an independence model
that satisfies the following axioms:
\begin{enumerate}[label=(G\arabic*),leftmargin=50pt]
\item\label{defn:gr:triv} (Triviality) $A\indep \emptyset\given C$;
\item\label{defn:gr:sym} (Symmetry) $A\indep B\given C \implies B\indep A\given C$;
\item\label{defn:gr:decomp} (Decomposition) $A\indep BD\given C \implies A\indep B\given C \AND A\indep D\given C$;
\item\label{defn:gr:wkun} (Weak union) $A\indep BD\given C \implies A\indep B\given CD$;
\item\label{defn:gr:contr} (Contraction) $A\indep B\given C \AND A\indep D\given BC \implies A\indep BD\given C$;
\item\label{defn:gr:int} (Intersection) $A\indep B\given CD \AND A\indep C\given BD \implies A\indep BC\given D$.
\end{enumerate}
The following additional axiom will be important in our discussion:
\begin{enumerate}[label=(G\arabic*),leftmargin=50pt]
\setcounter{enumi}{6}
\item\label{defn:gr:comp} (Composition) $A\indep B\given C \AND A\indep D\given C \implies A\indep BD\given C$.
\end{enumerate}
A graphoid that additionally satisfies \ref{defn:gr:comp} is called a \emph{compositional graphoid}. If only \ref{defn:gr:triv}-\ref{defn:gr:contr} are satisfied, then it is called a \emph{semi-graphoid}. A complete treatment of formal independence models and related concepts can be found in \cite{studeny2006probabilistic}.

Common structures that give rise to graphoids include probabilistic conditional independence, graphical separation, gaussoids, and \rev{partial correlation}. We will be interested in each of these structures in the sequel.

\medskip\noindent
\emph{Conditional independence.}
Let $\prob$ be a probability measure over the random vector $X=(X_{j})_{j\in\ground}$ and let $A,B,C$ be disjoint subsets of $V$. We write $X_{A}\indepprob X_{B}\given X_{C}$ to indicate that $X_{A}$ is conditionally independent of $X_{B}$ given $X_{C}$ under the probability measure $\prob$.
Denote the collection of all such conditional independence relations by $\ind(\prob)$. Then it is known that $\ind(\prob)$ is always a semi-graphoid, and if $\prob$ is a (strictly) positive measure then $\ind(\prob)$ is a graphoid.\footnote{A probability measure is called strictly positive if it has a density $q$ with respect to a product measure such that $q>0$; see e.g. Proposition~3.1 in \cite{lauritzen1996}.}
These (semi-)graphoids are called \emph{probabilistic (semi-)graphoids}\rev{, or simply \emph{probabilistic} for short}.
In many cases (e.g. Gaussian or symmetric binary), $\ind(\prob)$ is compositional, however, in general $\ind(\prob)$ need not be compositional \citep{studeny2006probabilistic,sadeghi2017faithfulness}. Each of the structures below furnishes examples of compositional graphoids.

\medskip\noindent
\emph{Graph separation.}
Let $G=(V,E)$ be an undirected graph (UG) and write $A\indepgr B\given C$ if and only if $A$ and $B$ are separated by $C$.
Denote the collection of all such separation statements by $\ind(G)$. Then $\ind(G)$ is a compositional graphoid, called the \emph{separation graphoid}. Similarly, separation graphoids may be defined for more general graphs including directed acyclic graphs (DAGs, via $d$-separation), chain graphs, mixed graphs, etc. \citep[cf.][]{lauritzen2018unifying}. In each case, the resulting separation graphoid is compositional \citep[e.g. Theorem 1 in][]{lauritzen2018unifying}. This has the important consequence that all results presented in the sequel apply to these general classes of graphical models.
If $\ind(\prob)=\ind(\gr)$, $\prob$ is called \emph{faithful} to $\gr$, and $\gr$ a \emph{perfect map} of $\prob$. In such cases, $\ind(\prob)$ is compositional.

\medskip\noindent
\emph{Gaussoids.}
Gaussoids are an abstraction of Gaussian independence models, introduced in \cite{lnenicka2007}. A gaussoid is a compositional graphoid that satisfies the additional axiom:
\begin{enumerate}[label=(G\arabic*),leftmargin=50pt]
\setcounter{enumi}{7}
\item\label{defn:gr:wktr} (Weak transitivity) $i\indep j\given C$ and $i\indep j\given Ck\implies i\indep k \given C$ or $j\indep k \given C$.
\end{enumerate}
Every regular Gaussian distribution $\normalN(0,\Sigma)$ gives rise to a gaussoid. However, not all gaussoids are representable by a Gaussian CI structure \citep{lnenicka2007, boege2019construction}, providing examples of compositional graphoids that are not necessarily (Gaussian) probabilistic or graphical.

\medskip\noindent
\emph{Partial correlation.}
\rev{Partial correlation, or second-order independence, also gives rise to a compositional graphoid (in fact, a gaussoid). This concept can be extended to general Hilbert spaces, \aarev{leading to what we refer to as the \emph{projection graphoid}}, discussed in detail in Section~\ref{sec:proj}.}

\medskip
\rev{It is a basic fact that not all graphoids are probabilistic, not all probabilistic graphoids are graphical, and not all compositional graphoids are graphical \citep[see e.g.][Section~3.6]{studeny2006probabilistic}. This motivates the need for a framework, like the one developed here, that can handle such general structures.}

\subsection{Examples}
\label{sec:ex}

To further motivate our results and illustrate the scope of compositional graphoids, we provide several explicit examples.

\begin{figure}[t]
    \centering
    \begin{tabular}{ccc}
        \includegraphics[width=2in]{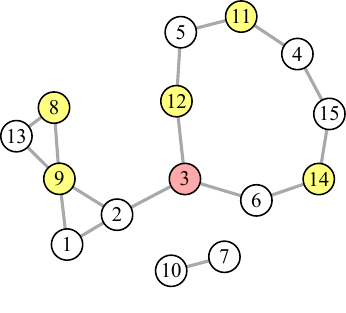}
    \end{tabular}
    \caption{\rev{Graph for Examples~\ref{exa:ug} and \ref{exa:bdy}-\ref{exa:perfect:cont} with $j=3$ and $S=\{8,9,11,12,14\}$ highlighted.}} 
\label{fig:exa:gr}
\end{figure}

\begin{example}[Graphical models]
\label{exa:ug}
Consider the graph $\gr$ shown in Figure~\ref{fig:exa:gr} with $d = 15$ nodes. We can associate to $\gr$ its separation graphoid $\graphoid(\gr)$. If $\prob$ is a Gaussian $\normalN(0,\Sigma)$ faithful to $\gr$ \citep[such a Gaussian always exists; see][]{lnenicka2007,amini2022perfectness}, then $\ind(\prob) = \ind(\gr)$ is a compositional graphoid. This graph will serve as a recurring example. \rev{More generally, our framework applies to other graphical model types (DAGs, chain graphs, etc.) as their separation properties also yield compositional graphoids.}
\end{example}
\noindent
The previous example can be substantially generalized:
\begin{lemma}
\label{lem:ex:cg}
If $\nu$ belongs to either of the following families of distributions, then $\ind(\nu)$ (the CI structure of $\nu$) gives rise to a compositional graphoid:
\begin{enumerate}
\item Distributions faithful to some (not-necessarily-parametric) graphical model $G$ (e.g., UG, DAG, chain graph, mixed graph, etc.). 
\item Regular Gaussian distributions, even those that are not faithful to any graph.
\end{enumerate}
\end{lemma}
\noindent
The first part follows from \citet[Theorem 1]{lauritzen2018unifying}. The second is a known property of Gaussian distributions \citep[e.g., Corollary~2.4 in][]{studeny2006probabilistic}.
Other examples of \update{compositional graphoids} include symmetric binary distributions, max-linear Bayesian networks, discrete determinantal point processes, and MTP2 distributions. Our framework applies to all such families without requiring \emph{a priori} knowledge of the specific model type or structure.

The following examples further motivate our work by illustrating the existence of non-graphical compositional graphoids. This demonstrates that, in addition to (a) subsuming existing graphical approaches and (b) removing the need to fix a particular graph type in advance, the neighbourhood lattice can model structures that \emph{cannot} be represented by a graph.

\begin{example}[A compositional graphoid that is not probabilistically representable]
\label{exa:studeny}
\cite{studeny1989multiinformation} provided an example of a compositional graphoid on $\ground=\{a,b,c,d\}$ defined by relations $a\indep b\given\{c,d\}$, $c\indep d\given a$, $c\indep d\given b$, $a\indep b$, their symmetric versions, and trivial independencies, which cannot be represented by any probability distribution (cf. Proposition~5 therein).
\end{example}

\begin{example}[Gaussoids that are not UGs]
\label{exa:gaussoid}
\cite{lnenicka2007} constructed gaussoids that are not representable by any Gaussian distribution. If such a gaussoid $\graphoid$ were representable by an undirected graph (UG) $\gr$, then by \citet[Corollary~3]{lnenicka2007}, there would exist a Gaussian $\prob$ such that $\graphoid(\prob)=\graphoid(\gr)=\graphoid$, a contradiction. Thus, gaussoids can furnish examples of non-graphical, compositional graphoids.
\end{example}

Other rich sources of examples are graphoids representable by one type of graphical structure but not others. While technically graphical, the practical challenge of knowing in advance which graph type to use motivates a more general approach.

\begin{figure}[t]
\centering
\begin{subfigure}[b]{0.4\textwidth}
\begin{tikzpicture}[node distance={15mm},
                    latent/.style = {draw, circle, orange},
                    obs/.style = {draw, circle, black},
                    obs2/.style = {draw, circle, red},
                    ]
\node[obs2] (X3) {$7$};
\node[obs2] (X1) [above left of=X3] {$5$};
\node[obs2] (X2) [above right of=X3] {$6$};

\node[obs] (X6) [left of=X1] {$3$};
\node[obs] (X4) [above left of=X6]{$1$};
\node[obs] (X5) [below left of=X4] {$2$};
\node[obs] (X7) [below left of=X6] {$4$};

\draw[-, black] (X4) -- (X5);
\draw[-, black] (X4) -- (X6);
\draw[-, black] (X5) -- (X7);
\draw[-, black] (X6) -- (X7);

\draw[->, red] (X1) -- (X3);
\draw[->, red] (X2) -- (X3);
\end{tikzpicture}
\caption{Chain graph in Example~\ref{exa:chain}.}
\label{fig:ex:chain}
\end{subfigure}
\qquad
\begin{subfigure}[b]{0.4\textwidth}
\quad\quad
\begin{tikzpicture}[node distance={15mm},
                    latent/.style = {draw, circle, orange},
                    obs/.style = {draw, circle, black}
                    ]
\node[latent] (H) {$H$};
\node[obs] (X2) [below left of=H] {$2$};
\node[obs] (X3) [below right of=H] {$3$};
\node[obs] (X1) [above left of=X2] {$1$};
\node[obs] (X4) [above right of=X3] {$4$};

\draw[->, black] (H) -- (X2);
\draw[->, black] (H) -- (X3);
\draw[->, black] (X1) -- (X2);
\draw[->, black] (X4) -- (X3);
\end{tikzpicture}
\caption{Hidden variable model in Example~\ref{exa:mixed}.}
\label{fig:ex:mixed}
\end{subfigure}
\caption{Examples of CI structures motivating general frameworks. }
\label{fig:ex}
\end{figure}

\begin{example}[A chain graph that is not a DAG or UG]
\label{exa:chain}
The chain graph in Figure~\ref{fig:ex:chain} (a $v$-structure combined with an undirected diamond) has a CI structure that cannot be represented by any single UG or DAG. \rev{Since Gaussian models of this form exist \citep{lnenicka2007}, this example is typical, not pathological} (see also Example~\ref{ex:unfaithful} for more).
\end{example}

\begin{example}[A mixed graph that is not a DAG or UG]
\label{exa:mixed}
Marginalizing a hidden variable $H$ from the DAG in Figure~\ref{fig:ex:mixed} results in a distribution $Q_H$ over $(X_1,X_2,X_3,X_4)$ whose CI structure $\graphoid(Q_H)$ cannot be perfectly mapped by any DAG or UG, though it can be represented by a mixed graph.
\end{example}

\section{Neighbourhood lattices and conditional independence}
\label{sec:lat}

Having introduced the intuitive concept of partitioning subsets based on shared explanatory properties (Section~\ref{sec:intro:intuitive_lattice}) and established the formal framework of compositional graphoids (Section~\ref{sec:bg}), we now rigorously define the neighbourhood lattice and explore its fundamental connection to conditional independence. This section lays the groundwork for understanding how these lattices provide a complete representation of the CI structure.
Throughout this section, we assume given a random vector $\rv$ with fixed index set $\ground$ whose CI relations induce a compositional graphoid $\graphoid$ over $\ground$.
In deliberate analogy with the usual convention for graphical models, we write $j\in\ground\iff\rv_{j}\in X$,
which will henceforth be referred to as \emph{nodes}. Everything in the sequel applies to general compositional graphoids over any finite ground set $\ground$.

\subsection{Relative Markov boundaries and blankets}
\label{sec:mb}

A cornerstone of understanding local dependencies is the concept of a Markov boundary.
The Markov boundary of an element $j \in \ground$ is typically defined as the smallest subset $U \subseteq \ground_{-j}$ such that $j\indep \ground_{-j}\setminus U\given U$~\citep{pearl1985graphoids,pearl1988}.
We generalize this by considering an arbitrary subset $S \subseteq \ground_{-j}$ as the current scope of interest.
The \emph{relative Markov boundary} of $j$ in (or with respect to) $S$ is then
\begin{align}
    \label{eq:def:bdy}
    \latmin(j;S)
    &= \bigcap\big\{ U\subseteq S : \; j \indep S\setminus U\given U\big\}.
\end{align}
In the sequel, we will frequently drop the qualifier ``relative'' when the context is clear, and use ``Markov boundary'' or simply ``boundary'' for brevity.
The boundary $\latmin(j;S)$ is well-defined for any $S$ over any graphoid, due to the intersection axiom \ref{defn:gr:int}. This follows more generally from the properties of the set of \emph{relative Markov blankets}.
For any subset $S\subseteq\ground_{-j}$, let
\begin{align}
\nhbdlatred_{j}(S)
= \{ U\subseteq S : \; j\indep S\setminus U\given U\}
\end{align}
be the collection of relative Markov blankets of $X_{j}$ with respect to $S$. The following lemma, whose proof is in Appendix~\ref{sec:proofs:lem:nhbdlatred:lat}, shows that $\nhbdlatred_{j}(S)$ has a desirable structure.

\begin{lemma}
\label{lem:nhbdlatred:lat}
Let $\graphoid$ be a graphoid over $\ground$, $j\in\ground$, and $S\subseteq\ground_{-j}$.
Then $\nhbdlatred_{j}(S)$ is an interval in $\Bf(\ground_{-j})$, and hence in particular a convex lattice.
\end{lemma}

\noindent
It follows that $\nhbdlatred_{j}(S)$ has a well-defined infimum, which is precisely $\latmin(j;S) = \inf\nhbdlatred_{j}(S)$, and a supremum, which is $S$ itself. Thus, $\nhbdlatred_{j}(S) = [\latmin(j;S), S]$.

\begin{example}
\label{exa:bdy}
Consider the graph $\gr$ from Example~\ref{exa:ug} (Figure~\ref{fig:exa:gr}). For $j = 3$, we can compute various relative Markov boundaries in the graphoid $\ind(\gr)$:
\begin{align*}
S_{1} = V_{-j} &\quad:\quad \latmin(j;S_{1}) = \{2,6,12\}, \\
S_{2} = \{2,6,12\} &\quad:\quad \latmin(j;S_{2}) = \{2,6,12\}, \\
S_{3} = \{4,6,12,14\} &\quad:\quad \latmin(j;S_{3}) = \{6,12\}, \\
S_{4} = \{7,10\} &\quad:\quad \latmin(j;S_{4}) = \emptyset, \\
S_5 = \{8, 9, 11, 12, 14\} &\quad:\quad \latmin(j;S_{5}) = \{9,12,14\}.
\end{align*}
In each case, $\latmin(j;S)$ is the smallest set $m \subseteq S$ separating $j$ from $S \setminus m$ in the graph $\gr$.
\end{example}

\subsection{The neighbourhood lattice}
\label{sec:nhbdlat}

The concept of a relative Markov boundary $\latmin(j;S)$ captures the essential variables within $S$ needed to shield $j$ from the rest of $S$. We now define a collection of sets $U$ that all share the \emph{same} relative Markov boundary with respect to $j$. This collection forms the core of our construction.

\begin{definition}
\label{defn:lattice}
Given a compositional graphoid $\graphoid$, a node $j\in\ground$, and $S\subseteq V_{-j}$,
define the \emph{neighbourhood lattice} of $j$ relative to $S$ as
\begin{align}
\label{eq:def:lattice}
\nhbdlat_{j}(S)
&= \{ U\subseteq\ground_{-j} :\;  \latmin(j;U) = \latmin(j;S) \}.
\end{align}
In words, $\nhbdlat_{j}(S)$ is the collection of all subsets $U \subseteq \ground_{-j}$ for which the relative Markov boundary of $j$ within $U$ is identical to the relative Markov boundary of $j$ within $S$.
When we wish to emphasize the underlying graphoid $\graphoid$, we will write $\nhbdlat_{j}(S;\graphoid)$.
\end{definition}

It follows from Lemma~\ref{lem:nhbdlatred:lat} (specifically, the well-definedness of $\latmin(j;S)$ and $\latmin(j;U)$) that $\nhbdlat_{j}(S)$ is well-defined.
Our first main result establishes the crucial properties of this collection, most notably that it indeed forms a lattice.

\begin{theorem}
\label{thm:nhbdlat}
Let $\graphoid$ be a compositional graphoid over $\ground$ and $j\in\ground$, $S\subseteq\ground_{-j}$ be fixed. Then:
\begin{enumerate}[label=(\alph*)]
\item\label{thm:nhbdlat:lat} $\nhbdlat_{j}(S)$ is an interval in $\Bf(\ground_{-j})$, and in particular a convex sublattice of it.
\item\label{thm:nhbdlat:bdy} $\inf\nhbdlat_{j}(S)=\latmin(j;S)$.
\item\label{thm:nhbdlat:sup} $\sup\nhbdlat_{j}(S) = S \cup E$, where $E=\{k\in\ground_{-j} \setminus S 
:\; k\indep j\given \latmin(j;S)\}$.
\item\label{thm:nhbdlat:decomp} The distinct neighbourhood lattices $\nhbdlat_{j}(S')$ (where $S'$ varies over $2^{\ground_{-j}}$) partition the subset lattice $\Bf(\ground_{-j})$. More precisely, if $m_S = \latmin(j;S)$, then $\nhbdlat_j(S) = \nhbdlat_j(m_S)$, and
\begin{align}
\label{eq:thm:nhbdlat:decomp}
2^{\ground_{-j}} = \biguplus\big\{\nhbdlat_{j}(\latmin) :\;
\aarev{\text{$m$ is minimal for $j$}}
\big\}.
\end{align}
\end{enumerate}
\end{theorem}
\noindent
\aarev{The statement ``$m$ is 
minimal for $j$'' is a shorthand for ``$\latmin=\latmin(j;T)$ for some $T\subseteq\ground_{-j}$''.}
Since $\nhbdlat_{j}(S)$ is an interval (and thus a lattice), it has a well-defined infimum and supremum. Theorem~\ref{thm:nhbdlat}\ref{thm:nhbdlat:bdy} characterizes the infimum as $\latmin(j;S)$, the common relative Markov boundary defining the lattice. The supremum, which we denote by
\begin{align}
\label{eq:def:latmax}
\latmax(j;S) := \sup\nhbdlat_{j}(S),
\end{align}
also plays an important role. Theorem~\ref{thm:nhbdlat}\ref{thm:nhbdlat:sup} characterizes this supremum by starting with $m = \latmin(j;S)$ and adding all other nodes $k \in \ground_{-j}$ that are conditionally independent of $j$ given $m$. Thus, each neighbourhood lattice $\nhbdlat_j(S)$ is precisely the interval $[\latmin(j;S), \latmax(j;S)]$.

\begin{example}
    \label{exa:lattice}
    Continuing Example~\ref{exa:bdy} with $j=3$ (see Figure~\ref{fig:exa:gr}),
\begin{align*}
    S_{1} = \ground_{-j}: &\quad \latmax(j;S_{1}) = \ground_{-j}. \\
    S_{2} = \{2,6,12\}: &\quad \latmax(j;S_{2}) = \ground_{-j}. \\
    S_{3} = \{4,6,12,14\}: &\quad \latmax(j;S_{3}) = \{4,5,6,7,10,11,12,14,15\}. \\
    S_{4} = \{7,10\}: &\quad \latmax(j;S_{4}) = \{7,10\}. \\
    S_5 = \{8,9,11,12,14\}: &\quad \latmax(j;S_{5}) =  \{4,5,7,8,9,10,11,12,13,14,15\}.
    \end{align*}
    Consequently, for instance, $\nhbdlat_{j}(S_{1}) = [\{2,6,12\},\ground_{-j}]$ and \[\nhbdlat_{j}(S_{5}) = [\{9,12,14\}, \{4,5,7,8,9,10,11,12,13,14,15\}]. \]
    The suprema $\latmax(j;S)$ are determined according to Theorem~\ref{thm:nhbdlat}\ref{thm:nhbdlat:sup}.
    Note that $\nhbdlat_j(S_1) = \nhbdlat_j(S_2)$ since $\latmin(j;S_1) = \latmin(j;S_2)$, illustrating that distinct $S$ can define the same neighbourhood lattice.
    The lattice $\nhbdlat_j(S_1)$ contains $2^{|\latmax(j;S_1)| - |\latmin(j;S_1)|} = 2^{14-3} = 2^{11}$ distinct subsets of $\ground_{-j}$. For any of these $2^{11}$ sets $U$, $\latmin(j;U) = \{2,6,12\}$. Similarly, $|\nhbdlat_j(S_5)|=2^{11-3} = 2^8=256$.
    \end{example}

\begin{remark}
    \label{rem:noncomp}
    Inspection of the proof of Theorem~\ref{thm:nhbdlat} (Appendix~\ref{sec:proof:main}) shows that if the graphoid $\graphoid$ is not compositional (i.e., fails axiom \ref{defn:gr:comp}), $\nhbdlat_{j}(S)$ is still closed under set intersections and possesses convexity. Composition is primarily needed for closure under set unions. Thus, in any graphoid, $\nhbdlat_{j}(S)$ is at least a meet semi-sublattice of $\Bf(\ground_{-j})$.
\end{remark}

\subsection{Implications for conditional independence}
\label{sec:ci}

The structure of the neighbourhood lattices, as established by Theorem~\ref{thm:nhbdlat}, has direct and powerful implications for determining CI relations. Specifically, knowing the interval $[m,M]$ that defines a neighbourhood lattice $\nhbdlat_j(S)$ allows us to read off many CI statements involving node $j$.

Our first result relates elementary CI relations of the form $j\indep i\given C$ to the neighbourhood lattice structure.
\begin{corollary}
\label{cor:elem:ci}
Let $\graphoid$ be a compositional graphoid over $\ground,$ $C\subseteq\ground_{-ji}$, and $i,j\in\ground$ with $i \ne j$. The following are equivalent:
\begin{enumerate*}[label=(\alph*)]
\item\label{cor:gen:ci:elem} $j\indep i\given C$;
\item\label{cor:gen:ci:elem:b} There is a neighbourhood lattice $\nhbdlat_{j}(S_0)=[\latmin,\latmax]$ (for some $S_0$, where $\latmin = \latmin(j;S_0)$) such that $i\in \latmax \setminus \latmin$ and $C\in[\latmin, \latmax \setminus i]$;
\item\label{cor:gen:ci:elem:c} $C \in \nhbdlat_{j}(Ci)$;
\item\label{cor:gen:ci:elem:d} $\latmin(j;Ci) \subseteq C$.
\end{enumerate*}
\end{corollary}
\noindent The proof is an application of the definitions and Theorem~\ref{thm:nhbdlat}.
Condition (d) is particularly useful: $j$ is independent of $i$ given $C$ if and only if $i$ is not part of the Markov boundary of $j$ relative to the set $Ci$.

To check general CI relations $A\indep B\given C$, we leverage the decomposition \ref{defn:gr:decomp} and composition \ref{defn:gr:comp} axioms, which imply that $A\indep B\given C$ if and only if $a \indep b \given C$ for all $a \in A, b \in B$. This leads to:
\begin{corollary}
\label{cor:gen:ci}
Let $\graphoid$ be a compositional graphoid over $\ground,$ and $A,B,C$ disjoint subsets of $\ground.$
Then, $A\indep B\given C$ if and only if $C \in \bigcap_{a\in A}\bigcap_{b\in B} \nhbdlat_{a}(Cb)$, or equivalently,
\begin{align}
\label{eq:cor:gen:ci:gen}
\bigcup_{a\in A} \bigcup_{b\in B} \latmin(a;Cb) \subseteq C.
\end{align} 
\end{corollary}
\noindent
These corollaries demonstrate a fundamental link: the neighbourhood lattice structure \emph{encodes} the CI relations. Given the lattices, CI queries are reduced to checking set inclusions and memberships. For an elementary CI $j \indep i \given C$, \eqref{eq:cor:gen:ci:gen} simplifies to $\latmin(j;Ci) \subseteq C$, matching Corollary~\ref{cor:elem:ci}\ref{cor:gen:ci:elem:d}. Conversely, given a lattice $\nhbdlat_{j}(S)=[\latmin,\latmax]$, we immediately infer $j\indep i\given C$ for all $i\in \latmax \setminus \latmin$ and all $C \in[\latmin, \latmax \setminus i]$.

\begin{example}
Reading off the graph in Example~\ref{exa:ug} (Figure~\ref{fig:exa:gr}), we can see that
$X_{3}\indep X_{9}\given (X_{1},X_{2})$ and
$(X_{1},X_{2})\indep(X_{8},X_{13})\given X_{9}$.
To check these via the corresponding neighbourhood lattice properties, we apply the conditions from Corollaries~\ref{cor:elem:ci} and~\ref{cor:gen:ci}.
For the first relation ($j=3, i=9, C=\{1,2\}$), Corollary~\ref{cor:elem:ci}\ref{cor:gen:ci:elem:d} requires $\latmin(j;Ci) \subseteq C$. Indeed,
\begin{align*}
\latmin(3; \{1,2,9\}) = \{1,2\} \subseteq C = \{1,2\},
\end{align*}
confirming the independence.
For the second relation ($A=\{1,2\}, B=\{8,13\}, C=\{9\}$), Corollary~\ref{cor:gen:ci} requires $\bigcup_{a\in A} \bigcup_{b\in B} \latmin(a;Cb) \subseteq C$. We have:
\begin{align*}
\bigcup_{a\in \{1,2\}}\bigcup_{b\in \{8,13\}} \latmin(a;\{9,b\})
&= \{9\} \cup \{9\} \cup \{9\} \cup \{9\} \subseteq C = \{9\},
\end{align*}
confirming the independence.
The individual $\latmin(\cdot\,;\cdot)$ terms are readily determined from the graph structure (e.g., $\latmin(1;\{9,8\}) = \{9\}$ because $X_9$ separates $X_1$ from $X_8$).
\end{example}

\subsection{The neighbourhood lattice decomposition}
\label{sec:nhbdlat:decomp}

Theorem~\ref{thm:nhbdlat}\ref{thm:nhbdlat:decomp} states that the distinct neighbourhood lattices partition the power set $2^{\ground_{-j}}$. This partition is central to our approach.

\begin{definition}[\aarev{Decomposition and sparsity}]
\label{defn:lat:decomp}
The \aarev{(full)} \emph{neighbourhood lattice decomposition} $\nhbdlatd_{j}$ of $j\in\ground$ is the partition of $\Bf(\ground_{-j})$ given by Theorem~\ref{thm:nhbdlat}\ref{thm:nhbdlat:decomp}:
\begin{align}
\label{eq:def:nhbdlat:decomp}
\nhbdlatd_{j}
= \big\{\nhbdlat_{j}(\latmin) : \;\aarev{\text{$m$ is minimal for $j$}}
    \big\}.
\end{align}
The $t$-\emph{sparse lattice decomposition} of $j$, denoted $\nhbdlatd_{j}(t)$, is
\begin{align}
\label{eq:def:nhbdlat:decomp:sparse}
     \nhbdlatd_{j}(t) = \big\{\nhbdlat_{j}(\latmin) : \; \text{$m$ is 
minimal for $j$ with $| m | \le t$} \big\}.
\end{align}
We say that node $j \in V$ has a \emph{$t$-sparse CI structure} if all its 
minimal sets
have size at most $t$, i.e., $\nhbdlatd_{j} = \nhbdlatd_{j}(t)$.
\end{definition}
\noindent
\aarev{In this section, we focus on the full decomposition $\nhbdlatd_j$; the consequences of sparsity are explored in subsequent sections, particularly Section~\ref{sec:highdim}.}
When node $j$ is clear, we may write $\nhbdlatd$. By \eqref{eq:thm:nhbdlat:decomp}, $2^{\ground_{-j}}=\bigcup \nhbdlatd$.

\aarev{The neighbourhood lattice decomposition $\nhbdlatd_j$ provides a compact encoding of all CI statements involving node $j$.}
As per Corollary~\ref{cor:elem:ci}, for each lattice $\nhbdlat^{\ell}=[\latmin^{\ell},\latmax^{\ell}]$ in $\nhbdlatd_j$, we have
\begin{align}\label{eq:all:CIs:listed}
    j \indep i \given  C \quad \text{for all} \quad
    i \in \latmax^{\ell}\setminus\latmin^{\ell} \text{ and for all } C \in [\latmin^{\ell},\latmax^{\ell}\setminus i],
\end{align}
and consequently:
\begin{corollary}\label{cor:ci:count} Let $\nhbdlatd_j = \{[m^\ell, M^\ell]\}_{\ell=1}^k$. Then, the total number of distinct elementary CI statements involving $j$ is $\sum_{\ell=1}^{k} (|\latmax^{\ell}| - |\latmin^{\ell}|)\, 2^{|\latmax^{\ell}|- |\latmin^{\ell}| - 1}$.  \end{corollary}
This number can be significantly less than enumerating all $(d-1) 2^{d-2}$ potential CI statements involving $j$. Thus, the decomposition offers a parsimonious representation, alternative to graphs, especially when a faithful graphical model is unknown or non-existent.

\begin{table}[t]
    \centering
    \begin{tabular}{lrrrrrrrrrr|r}
        & & & & & & & & & & & Total \\
        \hline
        Number of sets covered & 4 & 8 & 16 & 32 & 64 & 128 & 256 & 512 & 1024 & 2048 & $16384$\\ 
        \hline
        Number of lattices & 112 &  40 &  60 &  34 &  38 &  19 &   8 &   5 &   2 &   1 & $319$\\ 
        \hline
    \end{tabular}
    \caption{The distribution of the number of sets covered by each lattice $\nhbdlat = [m,M]$, i.e.  $2^{|M|-|m|}$, for the setup of Example~\ref{exa:lattice:full}.}\label{tab:covered:dist}
\bigskip
    \begin{tabular}{rrrrrrr}
        \hline
        size of $m$ & 0 & 1 & 2 & 3 & 4 & 5 \\ 
        \hline
        number of lattices &   1 &  12 &  60 & 134 &  91 &  21 \\ 
        \hline
    \end{tabular}
    \caption{The distribution of the size of the minimal element of the lattice $\Tc = [m,M]$, i.e. set $m$, for the setup of Example~\ref{exa:lattice:full}.}\label{tab:mcount:dist}
\end{table}

\begin{example}
\label{exa:lattice:full}
Continuing Examples~\ref{exa:bdy} and~\ref{exa:lattice} (Figure~\ref{fig:exa:gr}), for node $j=3$, the full neighbourhood lattice decomposition $\nhbdlatd_3$ consists of $319$ distinct lattices. Tables~\ref{tab:covered:dist} and~\ref{tab:mcount:dist} provide statistics. For instance, Table~\ref{tab:covered:dist} shows 5 lattices each covering $512$ subsets of $\ground_{-3}$. Table~\ref{tab:mcount:dist} shows that the largest relative Markov boundary size found for $j=3$ is $5$ (occurring for $21$ lattices). \aarev{Thus, node $j=3$ has a $5$-sparse CI structure according to Definition~\ref{defn:lat:decomp}.}
\end{example}

\section{Estimation and computation}
\label{sec:comp}

\newcommand{\Scand}[0]{candidates\xspace}
\newcommand{\minimals}[0]{minimals\xspace}
\let\oldReturn\Return
\renewcommand{\Return}{\State\oldReturn}

The theoretical results of Section~\ref{sec:lat}, particularly the structure of neighbourhood lattices (Theorem~\ref{thm:nhbdlat}) and their role in characterizing conditional independence (Corollaries~\ref{cor:elem:ci}-\ref{cor:gen:ci}), motivate the development of efficient estimation procedures. This section addresses two key computational questions:
\begin{enumerate}[label=(Q\arabic*)]
\item How can we efficiently determine the specific neighbourhood lattice $\nhbdlat_{j}(S)$ for a given node $j$ and subset $S$?
\item How can we efficiently compute the entire neighbourhood lattice decomposition $\nhbdlatd_{j}$ for a node $j$?
\end{enumerate}
We present algorithms for these tasks. To clearly separate  the algorithmic logic from statistical estimation issues, throughout this section, we assume access to a ``CI oracle'' that can answer queries of the form ``$A\indep B\given C$?'' for the given compositional graphoid $\graphoid$. In Section~\ref{sec:highdim}, we will discuss how to replace this oracle with statistical tests based on finite samples.

\subsection{Estimating Markov boundaries}
\label{sec:mbcomp}

A foundational step in our algorithms is the computation of the relative Markov boundary $\latmin(j;S)$ (defined in \eqref{eq:def:bdy}).
\rev{Many algorithms exist for learning Markov boundaries from data. For the general compositional graphoid setting here, a suitable procedure is an adaptation of the Grow-Shrink (GS) algorithm \citep{margaritis1999bayesian}.}
Algorithm~\ref{alg:bdy:comp} (detailed in Appendix~\ref{app:gs}) outlines such a procedure, \textsc{computeMB}$(j, S, \graphoid)$, which correctly finds $\latmin(j;S)$.
Its correctness relies only on the graphoid axioms \ref{defn:gr:triv}-\ref{defn:gr:comp} \citep[cf. Theorem~8,][]{statnikov2013algorithms}.
\textsc{computeMB} requires at most $O(|S|^{2})$ CI queries. \rev{While we use GS for complexity analysis, our framework can incorporate other correct MB learning algorithms.}

\subsection{Estimating a neighbourhood lattice}
\label{sec:latcomp}

The neighbourhood lattice $\nhbdlat_j(S)$ is the interval $[\latmin(j;S), \latmax(j;S)]$. Once the infimum $\latmin(j;S)$ is known, Theorem~\ref{thm:nhbdlat}\ref{thm:nhbdlat:sup} provides a direct way to find its supremum $\latmax(j;S)$.
This suggests Algorithm~\ref{alg:lattice:comp} for computing $\nhbdlat_{j}(S)$.
First, it calls \textsc{computeMB} to find $m = \latmin(j;S)$. 
Then, it initilizes the supremum as $M \gets S$,
and for each $k \in \ground_{-j} \setminus S$, adds $k$ to $M$ if $j \indep k \given m$. The algorithm returns the interval $[m,M]$.

\begin{figure}[t!]
    \centering
    \hfill
    \begin{minipage}[c]{.48\textwidth}
        \begin{algorithm}[H] \begin{algorithmic}[1]
                \Function{computeLattice}{$j$, $S$, $\graphoid$}
                \State $\latmin \gets$ \textproc{computeMB}($j$, $S$, $\graphoid$)
                \State $M \gets S$
                \State $A \gets \ground \setminus (S\cup\{j\})$ \ForAll{$k \in A$}
                \If {$j \indep k\given \latmin$}
                $M \gets M \cup \{k\}$
                \EndIf
                \EndFor
                \Return{$[\latmin,M]$}
                \EndFunction
            \end{algorithmic}
            \caption{Compute the lattice $\nhbdlat_{j}(S)$ for a given $j$ and $S \subseteq \ground_{-j}$.\label{alg:lattice:comp}}
        \end{algorithm}
    \end{minipage}\hfill
    \begin{minipage}[c]{.48\textwidth}
        \begin{algorithm}[H]
            \begin{algorithmic}[1]

\State $\nhbdlat \gets \textproc{computeLattice}(j, \emptyset,\graphoid)$
                \State Initialize $\nhbdlatd_j \gets \{\nhbdlat\}$.
                \State $\mathcal C  \gets $ the collection of all subsets of $V_{-j}$ of size at most $t$.
\For{$s=1,\ldots,t$}
                \For{$|S|=s$}
\State $\nhbdlat \gets \textproc{computeLattice}(j, S,\graphoid)$. 
\State $\nhbdlatd_j \gets \nhbdlatd_j \cup \{\nhbdlat\}$.
\EndFor
                \EndFor
\end{algorithmic}
            \caption{Compute the sparse lattice decomposition $\nhbdlatd_j(t)$ for a given $j$ and $t$. 
\label{alg:lattice:decomp:sparse}}
        \end{algorithm}
    \end{minipage}
\end{figure}

The \textsc{computeLattice} function involves one call to \textsc{computeMB} ($O(|S|^2)$ queries) and then $O(d)$ additional CI queries to find the supremum.
Thus, we have:
\begin{proposition}\label{prop:lat:comp:complex}
    The neighbourhood lattice $\nhbdlat_{j}(S)$ can be computed with $O(d+|S|^{2})$ CI queries.
\end{proposition}
In particular, computing a single neighbourhood lattice is polynomial in $d$.

\subsection{Computing lattice decompositions}
\label{sec:compdecomp}

We now consider computing the entire neighbourhood lattice decomposition $\nhbdlatd_j$ (Definition~\ref{defn:lat:decomp}), or its $t$-sparse version $\nhbdlatd_j(t)$.

\paragraph{Full Decomposition.}
The full decomposition $\nhbdlatd_j$ consists of all distinct lattices $\nhbdlat_j(m)$, where $m$ ranges over all possible relative Markov boundaries for node $j$.
Algorithm~\ref{alg:lattice:decomp} (detailed in Appendix~\ref{sec:proof:full:decomp}) provides a method to compute this full decomposition. It iteratively finds a subset $S \subseteq \ground_{-j}$ not yet covered by any discovered lattice, computes $\nhbdlat_j(S)$ using \textsc{computeLattice}, adds this unique lattice to the collection, and repeats until all of $2^{\ground_{-j}}$ is covered. The efficiency of finding an uncovered set relies on the disjointness of the lattices (Lemma~\ref{lem:decide:set:system}).

\begin{theorem}\label{thm:poly:N}
    The full lattice decomposition $\nhbdlatd_j$ can be computed with $O(d^{3}k^{2})$ CI queries, where $k=|\nhbdlatd_j|$ is the total number of distinct lattices in $\nhbdlatd_j$.
\end{theorem}
\noindent
This result implies that if the number of distinct lattices $k$ is relatively small, the entire CI structure around node $j$ can be learned efficiently.

\paragraph{$t$-Sparse Decomposition.}
Computing the full decomposition might be too costly if $k$ is large. An alternative is to compute the $t$-\emph{sparse} lattice decomposition $\nhbdlatd_j(t)$, which includes only those lattices $\nhbdlat_j(m)$ where the minimal element $m$ has size $|m| \le t$.
Algorithm~\ref{alg:lattice:decomp:sparse} outlines this approach. It simply iterates through all subsets $S\subseteq V_{-j}$ of size at most $t$, computes $\nhbdlat_j(S)$ and add it to the running set $\nhbdlatd_j$.
The use of a set data structure for $\nhbdlatd_j$ ensures that only distinct lattices are stored.    This is similar to how the PC algorithm learns the skeleton of a (faithful) DAG model \citep{spirtes2000,kalisch07}. In our case, of course, the faithfulness assumption is not necessary, and is replaced instead by the weaker composition axiom \ref{defn:gr:comp}. 

The computational complexity of Algorithm~\ref{alg:lattice:decomp:sparse} is determined by iterating over all candidate subsets $S \subseteq \ground_{-j}$ with $|S|=s$ for $s=0, \ldots, t$, and calling \textsc{computeLattice} for each. Since there are $\binom{d-1}{s}$ such subsets of size $s$, and each call to \textsc{computeLattice} costs $O(d+s^2)$ CI queries (Proposition~\ref{prop:lat:comp:complex}), the total complexity is $O\left( \sum_{s=0}^{t} \binom{d-1}{s} (d+s^2) \right)$. For fixed $t \ll d$, this sum is dominated by terms where $s \approx t$, leading to an overall complexity of $O(d^t(d+t^2))$, which can be further simplified to $O(d^{t+1})$ if $d \gg t^2$. This complexity is independent of $k$ (the total number of lattices in the full decomposition) but is exponential in $t$. If $j$ has a $t_0$-sparse CI structure (i.e., all its relative MBs are of size $\le t_0$), then running Algorithm~\ref{alg:lattice:decomp:sparse} with $t=t_0$ will recover the full decomposition $\nhbdlatd_j$, and $k \le \sum_{s=0}^{t_0} \binom{d-1}{s} \approx O(d^{t_0})$; this  bound could be too conservative as the next example shows. For the interested reader, \texttt{R} code implementing the algorithms is available at~\url{https://github.com/aaamini/nblat/}.

\begin{remark}
    \rev{The presented Algorithm~\ref{alg:lattice:decomp:sparse} iterates through all subsets $S$ up to size $t$. More efficient implementations for computing the sparse decomposition $\nhbdlatd_j(t)$ are possible by more carefully selecting candidate sets. For instance:
    1) One can maintain a record of regions in $2^{\ground_{-j}}$ already covered by found lattices. When considering a new candidate set $S$, if it falls within an already computed lattice, it can be skipped.
    2) If candidate sets $S$ are chosen strategically (i.e., by always picking an uncovered set $S$ of minimal possible cardinality), then $S$ itself would be the infimum $\latmin(j;S)$ for the new lattice $\nhbdlat_j(S)$ to be discovered. In such a case, the call to \textsc{computeMB}$(j, S, \graphoid)$ within the subsequent \textsc{computeLattice}$(j, S, \graphoid)$ procedure (i.e., Line 2 of Algorithm~\ref{alg:lattice:comp}) becomes redundant and can be optimized by directly setting $\latmin \gets S$.
    These optimizations can significantly reduce the number of CI queries, especially if many distinct $S$ values lead to previously discovered minimal boundaries $m$.}
\end{remark}

\begin{example}
    \label{exa:lattice:comp}
    Continuing Example~\ref{exa:lattice:full} (Figure~\ref{fig:exa:gr}), node $j=3$ has a $5$-sparse CI structure ($t_0=5$). The full decomposition has $k=319$ lattices.
    The naive bound $k \le \sum_{s=0}^5 \binom{14}{s} = 3,\!473$ is conservative.
    Computing the full lattice decomposition via Algorithm~\ref{alg:lattice:decomp} (from Appendix) for $j=3$ involved 2,208 CI queries.
    This should be contrasted with the total number of true CI statements involving $j=3$,  as given by Corollary~\ref{cor:ci:count}, which is 62,592 (out of $14 \cdot 2^{13} = 114,688$ possibilities).
    Thus, to deduce all 62,592 CI relations for $j=3$, only 2,208 CI queries and storing 319 lattice intervals are required, a substantial saving.
    If one were to use Algorithm~\ref{alg:lattice:decomp:sparse} with $t =5$, it would also recover the full decomposition.
\end{example}

\section{High-dimensional consistency}
\label{sec:highdim}

Let us now turn to estimating a CI structure (via the neighbourhood lattice decomposition) from a finite sample using the algorithms introduced in the previous section. The idea will be to use the same algorithms, but with the CI oracle replaced with a finite-sample test based on thresholding the generalized partial correlation, defined below.
Given the general nonparametric setting we are considering, testing even a single CI relation is known to be nontrivial \citep{shah2020hardness,neykov2020minimax}. Our goal in this section is to go far beyond a single such relation, and instead infer \emph{every} CI relation in $\prob$, or \aarev{equivalently},
the neighbourhood lattice decomposition which allows any CI relation to be efficiently checked (cf. Corollaries~\ref{cor:elem:ci}-\ref{cor:gen:ci}). 
Since there are exponentially many such relations, it is nontrivial to control the errors over all such tests, and this is where the lattice construction proves valuable.
\aarev{If $\nu$ is Gaussian, the problem can be reduced to a collection of linear neighbourhood regressions (cf. Section~\ref{sec:nbhd:lattice}) and techniques from sparse linear models applied. Without Gaussianity, however, linear regression techniques are not directly applicable. We show how ideas developed so far, coupled with a measure of CI for three scalar variables that concentrates well in high-dimensions, allows us to obtain results  that parallel the Gaussian case, for a fairly general class of distributions.}
\rev{Moreover, as in sparse linear regression, our results will depend on a sparsity parameter $t$, thereby avoiding the well-known curse of dimensionality in nonparametric estimation.}

First, we recall some definitions. To characterize conditional independence, we consider the \emph{generalized covariance measure} introduced in~\cite{shah2020hardness}, \aarev{although our results extend to any other measure for which a similar exponential concentration inequality holds}. 
Consider a pair of random variables $X_1, X_2 \in \reals$ and random vector $Z \in \reals^p$. The generalized covariance measure of $X_1$ and $X_2$, given $Z$, is
\begin{align}
    \gcm_{X_1, X_2 \given Z} := \ex[ \,\cov(X_1,X_2 \given Z)\, ].
\end{align}
It is clear that if $|\gcm_{X_1, X_2 \given Z}|  > 0$, then $X_1 \not\indep  X_2 \given Z$, although the converse is not necessarily true. For $j=1,2$, let $\hat X_j = \ex[X_j \given Z]$ be the regression function for predicting $X_j$ given $Z$, and let $X_j^\res = X_j - \hat X_j$ be the corresponding residual. It is not hard to show that
\[
\gcm_{X_1, X_2 \given Z}  = \cov(X_1^\res, X_2^\res) = \ip{X_1^\res, X_2^\res}_{L^2},
\]
where $\ip{X_1^\res, X_2^\res}_{L^2} = \ex[X_1^\res X_2^\res]$ is the usual $L^2$ inner product of two random variables.
We can now define the corresponding generalized correlation as
\begin{align}\label{eq:gcorr:def}
    \gcorr_{X_1,X_2 \given Z} := \frac{ \ip{X_1^\res, X_2^\res}_{L^2}}{\norm{X_1^\res}_{L^2}\norm{X_2^\res}_{L^2}}.
\end{align}
We refer to this quantity as \emph{generalized partial correlation} (GPC) coefficient.
The sample version of~\eqref{eq:gcorr:def} can be defined by using nonparametric regression to estimate $\ex[X_j \given Z]$ and then replacing population expectations in~\eqref{eq:gcorr:def} with the empirical ones. More precisely, assume that we observe an i.i.d. sample $(x_{i1},x_{i2},z_i) \sim (X_1,X_2,Z), \; i\in[n]$, and consider  estimates of the regression functions:
\begin{align}\label{eq:fh:j}
    \fh_j = \textup{argmin}_{f \in \Fc} \frac1n \sum_{i=1}^n(x_{ij} - f(z_i))^2
\end{align}
for $j=1,2$, where $\Fc$ is some function class. An estimate of $\gcm_{X_1,X_2 \given Z}$ is
\begin{align}\label{eq:def:gcmh}
    \gcmh_{X_1,X_2 \given Z} := 
    \frac1n \sum_{i=1}^n (x_{i1} - \fh_1(z_i)) (x_{i2} - \fh_2(z_i)),
\end{align}
and the corresponding estimate of the generalized correlation is
\begin{align}
    \gcorrh_{X_1,X_2 \given Z} := \frac{    \gcmh_{X_1,X_2 \given Z}}{(\gcmh_{X_1,X_1 \given Z} \cdot   \gcmh_{X_2,X_2 \given Z})^{1/2}}.
\end{align}

\newcommand\BCG{\operatorname{BCG}}

We are now ready to define the class of distributions we consider:
\begin{definition}
   The class of distributions $\nu$ such that
   \begin{enumerate*}[(a)]\item $\prob$ is supported on $[-1,1]^d$ and
        \item $\graphoid(\nu)$ is a compositional graphoid,
   \end{enumerate*}
   is denoted as $\BCG(d)$, and referred to as  $d$-dimensional bounded compositional graphoids. 
\end{definition}

\begin{definition}\label{defn:GPC:dist}
    Let $\alpha, \sigma > 0$, $t \in \nats$ and $\Fc$ be a class of functions from $\reals^t$ to $\reals$.
     A random vector $\rv=(\rv_{1},\ldots,\rv_{d})$ 
     has a $\gcmodel_t(\alpha, \sigma, \Fc)$ distribution $\prob$ if for any $i, j \in [d]$, and $S \subseteq [d]$ with $|S| \le t$, the following hold:
   \begin{enumerate}[(a),resume] \itemsep=.5ex
\item either $X_i \indep_{\!\!\prob} X_j \given X_S$ or 
            $
            |\gcorr_{X_i,X_j \given X_S}| \ge \alpha 
            $,
       \label{cond:alpha:sep}
\item $f^*_{i \given S} := \ex[X_i \given X_S] \in \Fc$,
        \item $\norm{X_i - f^*_{i \given S}}_{L^2} \ge \sigma$.
    \end{enumerate}
\end{definition}
\noindent
Condition~\ref{cond:alpha:sep} will also be referred to as ``$\alpha$-separated by generalized partial correlation, of order $t$.'' 
For the remainder of this section, we assume that $\rv=(\rv_{1},\ldots,\rv_{d})$ has a distribution which is both $\BCG(d)$ and $\gcmodel_t(\alpha, \sigma, \Fc)$. Let us write 
\begin{align}\label{eq:sample:gpcs}
\rhoh_n(i,j | S) = \gcorrh_{X_1, X_2 \given X_S}
\end{align}
for the corresponding sample GPC based on an i.i.d. sample $\rv^i=(\rv^i_{1},\ldots,\rv^i_{d}), \; i \in[n]$ from the distribution. Similarly, we write $\rho(i,j | S) = \gcorr_{X_1, X_2 \given X_S}$ for the population version.

Let $\nhbdlatd_j(t)$ be the population-level sparse lattice decomposition of order $t$ (cf. Definition~\ref{defn:lat:decomp}). One can obtain $\nhbdlatd_j(t)$ by replacing every query of the form $i \indep j \given S$ in Algorithm~\ref{alg:lattice:decomp:sparse} by checking whether or not $\rho(i,j|S) = 0$. Using the sample GPCs~\eqref{eq:sample:gpcs}, given a threshold $\tau > 0$, we can instead replace each such query with the test $|\rhoh_n(i,j|S)| \le \tau$. We refer to this version as the \emph{data-driven Algorithm~\ref{alg:lattice:decomp:sparse}} and denote its output as $\widehat \nhbdlatd_j(t, \tau)$.

The following result establishes the sample complexity needed so that with a data-driven choice of the threshold $\tau\asymp\sqrt{\eps_{n}^2+t\log d/n}$ (see \eqref{eq:crit:ineq:infnorm} for the definition of $\eps_{n}$), $\widehat \nhbdlatd_j(t, \tau)$ is a consistent estimate of $\nhbdlatd_j(t)$. 
For a function class $\Fc$ of functions from $\reals^d$ to $\reals$, we write $N_\infty(\eps, \Fc)$ for its $\eps$-covering number in the sup (a.k.a. uniform) norm and let $\Fc - \Fc =\{ f - g: f, g \in \Fc\}$. We say that a function class $\bar \Fc$ is star-shaped if whenever $f \in \bar \Fc$, then $\alpha f \in \bar \Fc$ for all $\alpha \in [0,1]$.

\begin{theorem}\label{thm:highd:consist:nongauss}
    Fix $\alpha, \sigma > 0$, $t \in \nats$ and let $\Fc$ be some function class.
    Suppose that  $\widehat \nhbdlatd_j(t, \tau)$ is computed based on an i.i.d. sample of size $n$ from a  distribution on $\reals^d$ which is both $\BCG(d)$ and $\gcmodel_t(\alpha, \sigma, \Fc)$.
Let $\Fcb$ be a star-shaped function 
class  containing $\Fc - \Fc$, and let $\eps_n$ be a solution of 
    \begin{align}\label{eq:crit:ineq:infnorm}
        C \int_{\eps^2/4}^\eps \sqrt{\log N_\infty(r, \Fcb)} \, dr \le \sqrt n \eps^2
    \end{align}
    where $C > 0$ is an absolute constant. There are absolute constants $C_1, C_2, C_3, c >0$ such that if
    \begin{align}\label{eq:alpha:range}
        C_1 \Bigl( \eps_n^2 + (t+4)\frac{ \log d}{n}\Bigr) \;\le\; \alpha^2 \;\le\; C_2 \min\{1, \sigma^2\},
    \end{align} 
    then, for any $\tau \in [\alpha/4, 3\alpha/4]$,
    \begin{align*}
        \pr\big( \widehat \nhbdlatd_j(t,\tau) = \nhbdlatd_j(t) \big) \ge 1 - C_3\, e^{-c \,n \alpha^2}.
    \end{align*} 
\end{theorem}

\noindent
The proof of Theorem~\ref{thm:highd:consist:nongauss} 
can be found in Appendix~\ref{app:proof_of_theorem_ref_thm_highd_consist}.
The interval $[\alpha/4, 3\alpha/4]$ can be replaced with $[c_1 \alpha, c_2\alpha]$ for any pair of constants $c_1$ and $c_2$ such that $0 < c_1 < c_2 < 1$, by modifying the other constants in the statement of the theorem. 
Since $\gcmodel_t(\alpha, \sigma, \Fc) \subseteq \gcmodel_t(\alpha', \sigma, \Fc)$ whenever $\alpha' \le \alpha$, if $\sigma = \Omega(1)$, the upper bound in~\eqref{eq:alpha:range} can always be satisfied by choosing a smaller $\alpha$. We assume $\sigma = \Omega(1)$ in the sequel and ignore the upper bound.

\revv{To explore the consequences of Theorem~\ref{thm:highd:consist:nongauss}, below we consider two examples: The nonparametric H\"older class and parametric Gaussian models. These examples illustrate two important aspects of this result: The H\"older class shows how the curse of dimensionality in nonparametric estimation can be avoided, and the Gaussian model provides an example of an unbounded function class for which Theorem~\ref{thm:highd:consist:nongauss} nonetheless still holds with the usual parametric rate.}

\begin{example}[\revv{H\"{o}lder classes}]
Let us consider a concrete function class, 
namely, let $\Fc$ be the set of functions from $[-1,1]^t$ to $\reals$ that belong to the unit ball of the H\"{o}lder space of order $s$. This is a standard example of a nonparametric function class that arises in statistics. Then, we can take $\Fcb$ to be the ball of radius 2 in the same space. Then, by Corollary~3 of~\cite{nickl2007bracketing} (with $\beta$ chosen to be $ > s$, which is valid since the domain is bounded), we have $\log N_\infty(r, \Fcb) \lesssim  r^{-t/s}$. 

Assume for simplicity that $1.9s \ge  t$. 
Then, $\int_0^\eps (\log N_\infty(r, \Fcb))^{1/2} dr \lesssim  \frac{2s}{2s-t} \eps^{1-(t/2s)} \le 20 \eps^{1-(t/2s)}$, and we can take
$
\eps_n^2 \lesssim n^{-2s / (2s+ t)}.
$
It follows that as long as 
\begin{align}\label{eq:alpha:scaling:holder}
\alpha^2 \gtrsim \Bigl( n^{-2s / (2s+ t)} + \frac{t \log d}{n} \Bigr),
\end{align}
we can recover the  $t$-sparse lattice decomposition of $j$, with high probability. As a result, we are able to infer not only all CI statements of the form $(i,j,S)$ with $|S| \le t$, but also many CI statements with $|S|>t$ as well (e.g. consider any lattice $\nhbdlat_{j}(S)$ such that $|\latmin_{j}(S)|=t$ and $|\latmax_{j}(S)|>t$). If it happens that $j$ has a $t$-sparse CI structure, then~\eqref{eq:alpha:scaling:holder} is enough to recover ``all'' CI statements involving $j$. 
\end{example}

\aarev{As discussed earlier, the $L^2$-norm-squared  rate of estimating an $s$-smooth function in $\reals^d$ is $n^{-2s/(2s+d)}$~\citep{stone1982optimal, bickel2007local, yang2015minimax} which gets extremely slow in high dimensions. 
The significance of the result obtained in~\eqref{eq:alpha:scaling:holder}}
is that instead of paying the price of full-dimensional curse of dimensionality $n^{-2s / (2s+ d)}$, we are paying the potentially much lower price of $n^{-2s / (2s+ t)}$. In particular, if  $j$ has a $t$-sparse CI structure with $t = O(1)$, we can recover all CI statements involving $j$ in the high-dimensional regime $d \gg n$ by only paying a $\log d$ price in dimension. Without leveraging the lattice decomposition, recovering even a single CI statement in the high-dimensional regime is hopeless, for the general distribution we consider here. 
\rev{For a more nuanced discussion of variable selection in nonparametric models, see \cite{CD11a,giordano2020grid}.}

\begin{example}[\revv{Gaussian models}]
\aarev{Unsurprisingly, a simplified version of Theorem~\ref{thm:highd:consist:nongauss} holds in the Gaussian case as well. GPC in this case reduces to the partial correlation coefficient (PCC).  The same proof goes through with the concentration inequality for GPC (Eqn.~\eqref{eq:rho:concent} in Appendix~\ref{sec:proofs:main}) replaced by a similar result for PCC, available in~\cite{kalisch07}. In the Gaussian case, the function class $\Fc$ can be taken to be the set of linear functions, leading to the parametric rate $\eps^2_n \asymp \frac{\log n}n$ in~\eqref{eq:alpha:range}. That is, for a Gaussian distribution, $\alpha$-separated by PCC, it is enough to have
\[
\frac{\log n}{ n} + (t+1) \frac{\log d}{n} \;\lesssim\; \alpha^2 \;\lesssim\; \min\{1,\sigma^2\}.
\]
This suggests that the general gap condition~\eqref{eq:alpha:range} has an optimal dependence on $\eps^2_n$ (that is, one that can recover the Gaussian rate). For brevity, we have omitted a formal statement in the Gaussian case; it is straightforward to supply one, if needed, based on the above discussion.}
\end{example}

\section{Graphical interpretation} 
\label{sec:grint}

The neighbourhood lattice and its corresponding minimal and maximal elements have intuitive interpretations in a graphical model; i.e. when $\graphoid$ is the separation graphoid of a graph $\gr$. Throughout this section we consider a fixed undirected graph $\gr=(V,E)$ and its separation relations $\ind(\gr)$. 
Then:
\begin{enumerate}[label=(B\arabic*)]
\item\label{grint:latmin} The relative boundary $\latmin(j;S):=\Ss$ of a node $j$ is the smallest subset of $S$ that separates $j$ and $S \setminus \Ss$;
\item\label{grint:latmax} The maximal set $\latmax(j;S)$ is obtained by adding to $\Ss$ every node $k\in\ground_{-j}$ that is separated from $j$ by $\Ss$;
\item\label{grint:inc} A subset $T$ is in the neighbourhood lattice $\nhbdlat_{j}(S)$ if (and only if) $\Ss\subseteq T$ and every $k\in T$ is separated from $j$ by $\Ss$.
\end{enumerate}

\noindent
These interpretations are direct translations of 
earlier definitions and results
to the special case of graphical separation, and recover familiar notions from the literature on graphical models. 
For example, when $S=\ground_{-j}$, \ref{grint:latmin} implies that
the Markov boundary is just the set of neighbours to $j$ in $\gr$. More generally, $\latmin(j;S)$ can be interpreted as a generalization of the concept of neighbourhood to arbitrary sets $S\subseteq\ground_{-j}$, i.e. $\latmin(j;S)$ is the set of neighbours to $j$ in $S$, in the sense that $\latmin(j;S)$ blocks $j$ from every node in $S$.
Moreover, Corollary~\ref{cor:gen:ci} simply says that $A$ and $B$ are separated by $C$ in $\gr$ if and only if every vertex in $A$ is separated from $B$ by $C$, which is precisely the definition of separation in an undirected graph.

We can further characterize the minimal and maximal sets of the neighbourhood lattice of $j$ via connected components of the graph resulting from the removal of $j$:
\begin{theorem}\label{thm:conn:comp} 
\sloppypar{Suppose that removing node $j$ and the edges connected to it breaks $\gr$ into $K$ connected components given by the vertex subsets $G_1,G_2,\dots,G_K \subseteq \ground_{-j}$. Let $S_k = S \cap G_k$. Then,}
    \begin{itemize}
        \item[(a)] $\latmin(j;S) = \biguplus_k \latmin(j;S_k)$,
        \item[(b)] $\latmax(j;S) = \bigcup_k \latmax(j;S_{k}) = \biguplus_k M(j;S_k,G_k)$, 
    \end{itemize}
    where $M(j;S_k,G_k)$ is the largest element of
$\nhbdlat_j(S_k;G_k) :=  \{ T \subseteq G_k :\; \latmin(j;T) = \latmin(j;S_k)\}$.
\end{theorem}
\noindent
    The proof of Theorem~\ref{thm:conn:comp} can be found in Appendix~\ref{sec:proof:thm:conn:comp}.
    Note that $\nhbdlat_j(S_k;G_k)$ is the lattice restricted to the ground set $G_k \cup \{j\}$ (i.e. instead of $\ground$). The original lattice given in Definition~\ref{defn:nhbd:lattice} can be thought of as $\nhbdlat_{j}(S_k;\ground_{-j})$.
For illustrative purposes, some simple consequences of Theorem~\ref{thm:conn:comp} are as follows:
    \begin{enumerate}[label=(B\arabic*)]
        \setcounter{enumi}{3}
        \item\label{grint:empty} $\latmin(j;S) = \emptyset$ iff there is no path between $j$ and $S$. (Separation by the empty set.)
        \item\label{grint:decomp} If $\gr$ decomposes into two disjoint components, say $\gr_{1}$ and $\gr_{2}$, then $\latmax(j;S)$ contains $\gr
        _{2}$ for any $j\in \gr_{1}$ and any $S$.\end{enumerate}

\begin{figure}
    \centering
    \begin{tabular}{ccc}
\includegraphics[width=1.8in]{figs/pcg_colored.pdf} \quad\quad &
\includegraphics[width=2in]{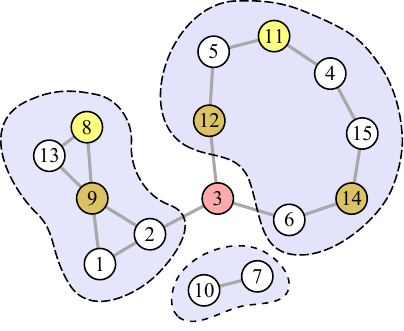}  \\
        (a) & (b) \end{tabular}
\caption{The graphical computation in Example~\ref{exa:perfect:cont} of the neighbourhood lattice $\nhbdlat_{j}(S)$ for $j=3$ and $S = \{9,8,11,12,14\}$. (a) Original graph $\gr$ with $j$ and $S$ specified with different colors. 
    From this figure, it is easy to see that $\Ss = \{9,12,14\}$ is the smallest subset of $S$ that separates $j$ from $S \setminus \Ss$, hence $m_j(S) = \Ss$. 
    (b)~Illustration of the three connected components which result after removing node $j=3$. As a result of Theorem~\ref{thm:conn:comp}, we can work in each component separately. For example, restricted to the left component $G_1 = \{1,2,8,9,13\}$, the minimal subset of $S \cap G_1 = \{8,9\}$ that separates $j=3$ from the rest of $S \cap G_1$ is $\Ss_1 = \{9\}$. Within $G_1$, $\Ss_1$ separates $\{8,13\}$ from $j=3$. Hence, the restricted lattice $\nhbdlat_{j}(S\cap G_1;G_1) = [\{9\},\{9,8,13\}]$.}
    \label{fig:pcg:graph:comp}
\end{figure}

\begin{example}\label{exa:perfect:cont}
Continuing the previous examples, let 
    $j=3$ and  $S =\{9,8,11,12,14\}$. 
It is clear from Figure~\ref{fig:pcg:graph:comp}(a) that the smallest subset $\Ss$ of $S$ separating $j$ from $S \setminus \Ss$ is $\{9,12,14\}$. Hence, $m_j(S) = \Ss$. Similarly, $E_j(\Ss) = \{4,5,7,8,10,11,13,15\}$ and $M_j(S) = \Ss \cup E_j(\Ss)$.
To verify Theorem~\ref{thm:conn:comp}, note that removing $j$ breaks $\gr$ into three connected components $G_1 = \{1,2,8,9,13\}$, $G_2 = \{7,10\}$ and $G_3 = \{4,5,6,11,12,14,15\}$, as illustrated in Figure~\ref{fig:pcg:graph:comp}(b). 
Then we can compute the restricted lattices $\Tc_j(S \cap G_k; G_k), k =1,2,3$. Using the notation $[m,M]$ to represent a lattice with minimum and maximum elements $m$ and $M$, respectively, the three lattices are: 
\begin{align*}
    [\{9\},\{8,9,13\}],\, 
    [\{12,14\},\{4,5,11,12,14,15\}],\, \text{ and } 
    [\emptyset, \{7,10\}]. 
    \end{align*}
    It is clear that the minimal and maximal elements of the original lattice are the disjoint union of the corresponding elements of these three lattices.
\end{example}

\begin{example}
\label{ex:mc}
Consider a Markov chain $X_{1}-X_{2}-\cdots-X_{d}$ with $\ground=[d]$. Define $j_{*}:=\sup\{\ell\in S : \ell < j \}$ and $j^{*}:=\inf\{\ell\in S : \ell > j \}$. Then we have for any $j$ and $S\subseteq\ground_{-j}$,
\begin{align*}
\latmin(j;S)
&= \begin{cases}
\{j_{*}\}, & \{\ell\in S : \ell > j \} = \emptyset \\
\{j^{*}\}, & \{\ell\in S : \ell < j \} = \emptyset \\
\{j^{*},j_{*}\}, & \text{otherwise},
\end{cases} \\
\latmax(j;S)
&= \begin{cases}
\{\ell\in\ground_{-j} : \ell \le j_{*}\}, & \{\ell\in S : \ell > j \} = \emptyset \\
\{\ell\in\ground_{-j} : \ell \ge j^{*}\}, & \{\ell\in S : \ell < j \} = \emptyset \\
\{\ell\in\ground_{-j} : \ell \ge j_{*}\}\cup \{\ell\in\ground_{-j} : \ell \le j^{*}\}, & \text{otherwise}.
\end{cases}
\end{align*}
To verify Theorem~\ref{thm:conn:comp}, there are three cases: (i) $j=1$, (ii) $j=d$, and (iii) $1<j<d$. In (i) and (ii), removing $j$ results in a single connected component, corresponding to the first two cases above. In (iii), removing $j$ leaves two connected components, corresponding to the third case.
\end{example}

\begin{example}
\label{ex:unfaithful}
Recall Example~\ref{exa:chain}, which describes a graphoid that is probabilistically representable but not representable as the separation graphoid of any undirected graph or directed acyclic graph.
We can explicitly construct an example by considering $\prob=\normalN(0,\Sigma)$ where $\Sigma$ is defined by: 
\begin{align*}
\Sigma^{-1} = \begin{pmatrix}
3 & 1 & 1 & 0 & 0 & 0 & 0 \\
1 & 3 & 0 & 1 & 0 & 0 & 0 \\
1 & 0 & 3 & 1 & 0 & 0 & 0 \\
0 & 1 & 1 & 3 & 0 & 0 & 0 \\
0 & 0 & 0 & 0 & 4 & 3 &-1 \\
0 & 0 & 0 & 0 & 3 & 5 &-3 \\
0 & 0 & 0 & 0 &-3 &-3 & 3 \\
\end{pmatrix}.
\end{align*}
Consider the undirected graph $G=(V,E)$ given by the nonzero entries of $\Sigma^{-1}$. It is easy to verify that $X_5\indep_{\!\!\!\prob} X_6$ but are not separated by the empty set (i.e. $(5,6)\in E$). That is, $(5,6,\emptyset)\in\ind(\prob)$ but $(5,6,\emptyset)\notin\ind(\gr)$.
Thus $\ind(\gr)\ne\ind(\prob)$.  
Indeed, as described in Example~\ref{exa:chain}, no undirected graph satisfies $\ind(\gr)=\ind(\prob)$.
Nonetheless, the lattice construction still applies since $\ind(\prob)$ is a gaussoid.
For example, suppose we wish to check the CI relation $X_5\indep_{\!\!\!\prob}X_6$, which can be done via \eqref{eq:cor:gen:ci:gen} with $j=5$, $i=6$, and $C=\emptyset$:
\begin{align*}
\latmin(j;Ci)=
\latmin(5;\{6\})=
\emptyset
\subseteq C = \emptyset.
\end{align*}
Since $\prob$ is a Gaussian, it is easy to verify that indeed $X_5\indep_{\!\!\!\prob}X_6$.
\end{example}

\section{Projection and regression interpretation}
\label{sec:proj}

An alternative perspective on the neighbourhood lattice arises in the context of regression models in statistical applications. In particular, the notion of \emph{neighbourhood regression} has played a prominent role in learning graphical models (e.g. \cite{meinshausen2006,ravikumar2010,yang2015}). This in turn bears a natural relation to orthogonal projections, which suggests a possible connection between the lattice of projections on a Hilbert space and the neighbourhood lattice. 
\rev{This connection gives rise to the concept of \emph{partial orthogonality}, closely related to partial regression and second-order independence, which we will define before illustrating its close connection with neighbourhood regression.}
We then explore these connections and illustrate the utility of the projection interpretation by providing an independent proof (based on abstract projection lattices) that partial orthogonality gives rise to a neighbourhood lattice.

\subsection{Partial orthogonality}
\label{sec:po}
Let $\Hil$ be a Hilbert space and let $\Prj(\Hil)$ denote its projection lattice, i.e., the set of (bounded) projection operators acting on $\Hil$. \rev{Section~\ref{sec:background:proj} contains more background on the projection lattice for interested readers.} We write $\csp(\Hilset)$ for the closed linear span of any subset $\Hilset$ of $\Hil$. Note that $\csp(\Hilset)$ can be viewed as a Hilbert space with the inner product inherited from $\Hil$. We write $\Prj(\Hilset) :=  \Prj(\csp(\Hilset))$ for the projection lattice of  $\csp(\Hilset)$.
    
Let $\Hilset$ be a finite subset of linearly independent vectors in $\Hil$.
For any $S \subseteq\Hilset$, let $P_S \in \Prj(\Hil) $ denote the projection onto $\csp(S)$
and $P_S^\perp := I - P_S$ denote the orthogonal complement of $P_S$. We note that $\{P_S :\; S \subseteq \Hilset\}$ is exactly $\Prj(\Hilset)$.
For any disjoint triple of subsets $A,B,C \subseteq\Hilset$, we say that $A$ is partially orthogonal to $B$ given $C$ if 
\begin{align}\label{eq:partial:orth}
    P_A P_C^\perp P_B = 0,
\end{align}
in which case we also write $P_A \perp P_B \mid P_C$ or $(A \perp B \mid C)$.
We can view~\eqref{eq:partial:orth} as a ternary relation among projection operators and say that ``$P_A$ is partially orthogonal to $P_B$ given $P_C$''.
\rev{When $\Hil$ is an $L^2$ space of square-integrable random variables on some probability space and $\Hilset$ is a finite set of random variables (i.e. a random vector), \eqref{eq:partial:orth} is equivalent to partial uncorrelatedness of $A$ and $B$ given $C$, also equivalent to second-order independence (Example~2.26 in \cite{lauritzen2020lectures}).}
\rev{This relation is also closely related to the \emph{orthogonal meet} found in prior work;} our choice of terminology is intended to reflect the well-known relationship with partial regression and partial correlation.
Equivalently, we can view partial orthogonality as a ternary relation on families of closed subspaces of $\Hil$;
e.g. as in \cite{dawid2001separoids}.
The resulting relation defines a gaussoid (see Appendix~\ref{sec:details:nbhd:lat} for a proof) which we call the \emph{projection graphoid} and denote by $\ind(\Prj(\Hilset))$. Since $\ind(\Prj(\Hilset))$  is a gaussoid, it is in particular a compositional graphoid over $\ground=\Hilset$.

\subsection{The regression lattice}
\label{sec:nbhd:lattice}

Let $X=(X_{j})_{j\in\ground}$ be a square-integrable random vector with covariance matrix $\Sigma=\cov(X)$ and joint distribution $\prob$. We view each $X_j$ as an element of the $L^2$ space of random variables. For any $S \subseteq \ground$, let $P_S$ denote the $L^2$ projection onto the span of $\gv_S = \{\gv_i,\; i\in S\}$.  For simplicity, we write $P_j = P_{\{j\}}$ to denote projection on the span of $\{\gv_j\}$, i.e., we drop the brackets for singleton sets.
With some abuse of notation, let us also view $X$ as a set $X = \{X_j, j \in \ground\}$ and recall the graphoid $\ind(\Prj(X))$ defined by the partial orthogonality relation~\eqref{eq:partial:orth}. Relative to $\ind(\Prj(X))$, we have Markov boundaries $m(j;S)$ and the neighbourhood lattices $\nhbdlat_j(S)$ defined via~\eqref{eq:def:bdy} and~\eqref{eq:def:lattice}, respectively. Throughout this section, in order to avoid ambiguity, we always explicitly reference the underlying graphoid, i.e. by writing $\nhbdlat_{j}(S;\prob) := \nhbdlat_{j}(S;\graphoid(\prob))$ for the probabilistic graphoid induced by conditional independence and $\nhbdlat_{j}(S;\Prj(X)) := \nhbdlat_{j}(S;\graphoid(\Prj(X)))$ for the projection graphoid induced by partial orthogonality (cf. \eqref{eq:partial:orth}). Note that these are in general distinct.

The following, which we call the \emph{regression lattice}, is closely related to the neighbourhood lattice $\nhbdlat_{j}(S;\Prj(X))$:

\begin{definition}[Regression lattice]
    \label{defn:nhbd:lattice}
    For any $S \subseteq \ground_{-j}$, define a collection of subsets by
    \begin{align}\label{eq:lattice:proj:def}
    \Tc_j(S) = \{ U \subseteq \ground_{-j}:\; P_U P_j = P_S P_j\}.
    \end{align}
\end{definition}
\noindent
In other words, for any $j$, $\nhbdset_j(S)$ is the collection of subsets $U\subseteq\ground_{-j}$ such that the projection of $\dvec_{j}$ onto $\{\dvec_i, \, i\in U\}$ is invariant. 
Note that $\Tc_j(S)$ is well-defined even when $\Sigma$ is rank-deficient.   

The name ``regression lattice'' may seem peculiar by inspection of the definition alone. That $\Tc_j(S)$ is indeed a lattice will be established in Theorem~\ref{thm:reg:lattice} below. Let us first explain our use of the adjective ``regression'':
Let $\beta_j(S)$ be the coefficient vector obtained after regressing $X_j$ onto $X_S$, that is
\begin{align}
\label{eq:def:nhbd:regression}
\coef_j(S) &:= \textup{argmin}_{\beta \,\in\, \R^d, \;\supp(\beta)\, \subseteq\, S}
\ex(\dvec_{j} - \beta^T\dvec)^{2},
\end{align}
where $\supp(a)=\{i : a_{i}\ne 0\}$ denotes the support of a vector $a$.
Equivalently, $\coef_j(S)^{T}X$ is the orthogonal projection of $X_{j}$ onto $X_{S}$ in the space of square-integrable random variables.
The coefficients in $\coef_j(S)$ are also known as the \emph{partial regression coefficients} for the variable $j$ regressed on the variables $S$. We will often refer to them simply as \emph{regression coefficients}.
Then, assuming $\Sigma \succ 0$, we have (cf.~Lemma~\ref{lem:equiv:lattice:defs} in Appendix~\ref{sec:details:nbhd:lat}): \begin{align}
\begin{split}
\nhbdset_j(S) &= \{U \subseteq\ground_{-j}: \; \nhbdcoef_j(U) = \nhbdcoef_j(S) \} \\
&= \big\{U \subseteq\ground_{-j}: \; \supp(\nhbdcoef_j(U)) = \supp(\nhbdcoef_j(S)) \big\}.
\label{eq:lattice:alt}
\end{split}
\end{align}
By definition, $ \supp(\beta_j(S)) \subseteq S$.

The next result shows that $\Tc_j(S)$ is in fact the same as the neighbourhood lattice induced by partial orthogonality with respect to the random vector $X$. As an immediate consequence, this also proves that $\Tc_j(S)$ is indeed a lattice.

\begin{theorem}
\label{thm:reg:lattice}
    Assume that $\Sigma := \cov(X) \succ 0$. Then, we have 
    \begin{align}
\nhbdset_j(S) = \nhbdlat_j(S;\Prj(X))
    \end{align}
    for all $j \in \ground$ and $S \subseteq \ground_{-j}$. \end{theorem}
\noindent
The proof of Theorem~\ref{thm:reg:lattice} can be found in Appendix~\ref{sec:proj:proofs}.
Combined with \eqref{eq:lattice:alt}, Theorem~\ref{thm:reg:lattice} explains the name \emph{regression lattice}. In the remainder of this section, we explore the properties of this lattice and its connection to the previously defined neighbourhood lattice.

\begin{remark}
By definition, the regression lattice $\Tc_j(S)$ depends only on the random vector $X$, and is not defined relative to any underlying graphoid. Contrast this with the definition of the neighbourhood lattice, which explicitly requires an underlying compositional graphoid $\graphoid$ that may not be generated by the distribution of some random vector $X$. In particular, we stress that $\Tc_j(S)$ is \emph{not} the same as $\nhbdlat_{j}(S;\prob)$, \update{where $\prob$ is the distribution of $X$.} See also the discussion after Lemma~\ref{lem:lattice:compare} below.
\end{remark}

\subsection{Comparison of lattices}
\label{sec:comparison}

We have now defined two objects: The \emph{neighbourhood} lattice and the \emph{regression} lattice. 
Theorem~\ref{thm:reg:lattice} showed that the regression lattice is in fact the same as the neighbourhood lattice $\nhbdlat_j(S;\Prj(X))$ that arises from partial orthogonality, as defined in \eqref{eq:partial:orth}. In this section, we explore the relationship between the regression lattice and the neighbourhood lattice $\nhbdlat_j(S;\prob)$ that arises from conditional independence.

Consider the setting in which we have a random vector $X\in L^{2}$, so that both $\nhbdlat_{j}(S;\prob)$ and $\Tc_{j}(S)$ are well-defined lattices. It turns out that $\nhbdlat_{j}(S;\prob)$ is always a subset of $\Tc_{j}(S)$, and equality holds if $X$ is Gaussian.

\begin{lemma}
\label{lem:lattice:compare}
Let $X=(X_{j})_{j\in\ground}\in L^{2}$ be a random vector and \update{$\prob$ be the distribution of $X$}. Then:
\begin{enumerate}[label=$(\alph*)$]
\item\label{} Without additional assumptions, $\nhbdlat_j(S;\prob) \subseteq \Tc_j(S)$.
\item\label{} If $X\sim\normalN(0,\Sigma)$ for some $\Sigma\posdef0$, then $\nhbdlat_j(S;\prob) = \Tc_j(S)$.
\end{enumerate}
\end{lemma}
\noindent
Thus, when $X\sim\normalN(0,\Sigma):=\prob$, we have
\begin{align}
\label{eq:gaussian:lat}
\nhbdlat_{j}(S;\prob)
=\mathcal{T}_{j}(S)
=\nhbdlat_{j}(S;\Prj(X)).
\end{align}
This is not true in general, however, since partial orthogonality is weaker than conditional independence. \revv{Thus, Lemma~\ref{lem:lattice:compare} simply captures the well-known relationships between independence, uncorrelatedness, and Gaussianity.}

\subsection{Projection lattice approach}
\label{sec:projlattice}

We now give an independent proof that $\Tc_j(S)$ is a convex lattice, using only the algebraic properties of projections. The proof reveals where the lattice structure of $\Tc_j(S)$ comes from: It is inherited from the natural lattice structure possessed by projections on a Hilbert space. We will assume that the reader is familiar with the basic theory of Hilbert spaces and their associated projection lattices; a detailed introduction to these topics can be found in \cite{Farah2010,Blackadar2006}. We provide a short review in Section~\ref{sec:background:proj} below.

\subsubsection{Background on projection lattices}\label{sec:background:proj}
For a (separable) Hilbert space $\Hil$, let $B(\Hil)$ be the space of bounded linear operators on $\Hil$. For an operator $P \in B(\Hil)$, let $\ran(P) := P \Hil := \{P x :\; x \in \Hil\}$ denote its range and $P^*$ its adjoint, defined via the relation $\ip{x,P^*y} = \ip{Px,y}$ for all $x,y \in \Hil$.
An operator $P \in B(\Hil)$ is an orthogonal projection if and only if it is self-adjoint and idempotent: $P^* = P = P^2$. The set of orthogonal projections in $B(\Hil)$ is denoted as $\Prj(\Hil)$. The range of any orthogonal projection is a closed linear subspace of $\Hil$. In fact, there is a bijection between $\Prj(\Hil)$ and closed linear subspaces of $\Hil$. The latter can be ordered by inclusion, which induces a natural (partial) order on $\Prj(\Hil)$ via the bijection. This order can be characterized as follows:
\begin{lemma}\label{lem:proj:order}
For $P,Q \in \Prj(\Hil)$, the following are equivalent:
    \begin{align*}
        (a) \; P = PQ, \quad (b)\; P = QP, \quad (c)\;\ran(P) \subseteq \ran(Q).
    \end{align*}
\end{lemma}
\noindent    When any of these conditions hold we write $P \le Q$.

The equivalence of (a) and (b) follows from self-adjointness of orthogonal projections. Note that 
$PQ$ in Lemma~\ref{lem:proj:order} is not necessarily a projection (unless $P$ and $Q$ commute). One can show that the above order turns $\Prj(\Hil)$ into a complete lattice. The meet and join of two elements $P,Q \in \Prj(\Hil)$ can be expressed as follows:
\begin{align*}
    P \wedge Q &=  \text{the projection onto $\ran(P) \cap \ran(Q)$} \\
    P \vee Q &= \text{the projection onto the closed linear span of $\ran(P) \cup \ran(Q)$.}
\end{align*}
We also let $P^\perp  := I - P$, the orthogonal complement of $P$. Note that $P \vee P^\perp = I$ (the identity operator) and $P \wedge P^\perp = \{0\}$. Also, $P \le Q$ iff $Q^\perp \le P^\perp$. See for example~\citet{Farah2010} (Section 5, p.~24) or \citet{Blackadar2006} (Section~II.3.2, p.~78) and the references therein.

\subsubsection{Abstract lattice theorem}
\label{sec:abstractlat}

We now show that $\Tc_j(S)$ is a lattice isomorphic to an interval of the subset lattice, by combining two abstract results. Consider a Hilbert space $\Hil$ and a collection of vectors $\mathcal X = \{x_j\}_{j \in J} \subseteq \Hil$, not necessarily finite. 
The next result and its proof (except the convexity assertion) are due to \citet{Tristan16}:
\begin{theorem}\label{thm:abs:proj}
    For an operator $A \in B(\Hil)$ and projection $P \in \Prj(\Hil)$, define
    \begin{align*}
        \mathcal Q(A, P) := \{Q \in \Prj(\Hil):\, P A = Q A\}.
    \end{align*}
    Then $\mathcal Q(A, P)$ is a  complete convex sublattice of $\Prj(\Hil)$.
\end{theorem} 
\noindent
The proof of Theorem~\ref{thm:abs:proj} can be found in Appendix~\ref{sec:proj:proofs}.

For any finite $S \subseteq \Xc$, let $P_S \in \Prj(\Hil)$ be the projection onto the (closed) linear span of $S$. Recall that $\Bf(\Xc)$ is the complete lattice of all subsets of $\Xc$ ordered by inclusion. We also write $\Bfin(\Xc)$ for the lattice of all finite subsets of $\Xc$. 
\begin{proposition}\label{prop:geom:lat:isom}
    Assume that $\Xc \subseteq \Hil$ is a linearly independent set. Then, 
    \[
    \Pi(\Xc) := \{P_S:\; S \; \text{is a finite subset of}\; \Xc\}
    \]
    is a sublattice of $\Prj(\Hil)$ isomorphic to  $\Bfin(\Xc)$.
\end{proposition}
\noindent
The proposition implies that for any $S, T \subseteq \Xc$, we have $P_S \wedge P_T = P_{S \cap T}$ and $P_S \vee P_T = P_{S \cup T}$. 
Combining these two results, and noting that $\Tc_j(S) = \Pi(\{X_i\}_{i \in \ground_{-j}}) \cap \mathcal Q(P_{X_j}, P_S)$,
we have the following corollary:
\begin{corollary}
$\Tc_j(S)$ is lattice-isomorphic to an interval of the subset lattice.
\end{corollary}

\begin{remark}
Proposition~\ref{prop:geom:lat:isom} and Theorem~\ref{thm:abs:proj} apply to an infinite-dimensional Hilbert space $\Hil$ and infinitely many variables $\{\gv_1,\gv_2,\dots\}$. They also hold when the variables are dependent (i.e., $\Sigma$ is rank-deficient) if we remain at the level of projections (i.e., not map projections onto sets of variables).

\end{remark}

\section{Discussion}\label{sec:discuss}

We have introduced the neighbourhood lattice for general compositional graphoids and its application to learning CI structures. Notably, the neighbourhood lattice obviates the need for a graphical representation of a distribution. We have shown that these lattices are efficiently computable, and have meaningful interpretations in special cases such as graphs and regression. \rev{We also established a high-dimensional consistency result that applies to general models, including various nonparametric models as well as more conventional Gaussian models. In both cases, the lattice structure significantly reduces the statistical complexity of the estimation problem.}

We conclude with some additional open questions. In Remark~\ref{rem:noncomp}, we observed that many of the properties carry over even for non-compositional graphoids, however, a complete study of this more general case remains open. It would also be interesting to derive sharper upper bounds on the number of lattices in the lattice decomposition (Definition~\ref{defn:lat:decomp}), as this has important computational implications (Section~\ref{sec:comp}). It would also be interesting to explore consequences of the neighbourhood lattice for learning graphical models, e.g. as in \cite{aragam2019globally}. In particular, since our results apply to \emph{any} graphical model, this opens the door for learning more general (e.g. chain, mixed, etc.) graphical models without requiring faithfulness.

\paragraph{Acknowledgements}
{AA was supported by NSF DMS-1945667. BA was supported by NSF IIS-2453378, IIS-1956330, and NIH R01GM140467. QZ was partially supported by NSF DMS-2305631. Part of this research was performed while BA was visiting the Institute for Mathematical and Statistical Innovation (IMSI), which is supported by the National Science Foundation (Grant No. DMS-2425650).}

\newpage

\appendix
\section{Alternative characterizations}
\label{sec:alt:char}

The following result lists several different equivalent characterizations of the neighbourhood lattice:
\begin{theorem}
    \label{thm:nhbdlat:equiv}
    The following definitions of $\nhbdlat_{j}(S)$ are equivalent:
    \begin{enumerate}[label=(A\arabic*)]
        \item\label{thm:nhbdlat:equiv:1} $\nhbdlat_{j}(S) = \{ U\subseteq\ground_{-j} : \latmin(j;U) = \latmin(j;S) \} = \{ U\subseteq\ground_{-j} : \inf \mathfrak{t}_{j}(U) = \inf \mathfrak{t}_{j}(S)\}$;
\item\label{thm:nhbdlat:equiv:2} $\nhbdlat_{j}(S) = \bigcup \{ \mathfrak{t}_{j}(U) : \inf \mathfrak{t}_{j}(U) = \inf \mathfrak{t}_{j}(S)\}$, i.e. $\nhbdlat_{j}(S)$ is maximal amongst the lattices $\mathfrak{t}_{j}(S)$ with the same minimal element;
\item\label{thm:nhbdlat:equiv:3} $\nhbdlat_{j}(S) = \{U\subseteq S : j\indep S-U\given U\} \cup \{U\supset S : j\indep U-S\given S\}$;
        \item\label{thm:nhbdlat:equiv:4} $\nhbdlat_{j}(S) = \mathfrak{t}_{j}(S) \cup \{ A\cup B : A\cap B=\emptyset,\, A\in\mathfrak{t}_{j}(S),\, j\indep B\given A\}$;
        \item\label{thm:nhbdlat:equiv:5} $\nhbdlat_{j}(S) = [m,M]$, where $m\subseteq\ground$ is the smallest subset such that $j\indep S-m\given m$ and $M\subseteq\ground$ is the largest subset such that $j\indep M-S\given S$.
        \item\label{thm:nhbdlat:equiv:6} $\nhbdlat_{j}(S) = \{T\subseteq\ground_{-j} : j\indep T-\latmin(j;S)\given \latmin(j;S)\}$.
\end{enumerate}
    Furthermore, we have the following property: 
\begin{enumerate}[label=(A\arabic*)]
        \setcounter{enumi}{6}
        \item\label{thm:nhbdlat:equiv:7} $\nhbdlat_{j}(S)=\nhbdlat_{j}(T)$ for all $T\in\nhbdlat_{j}(S)=[m,M]$ as in \ref{thm:nhbdlat:equiv:5}.
\end{enumerate}
\end{theorem}
\noindent
The proof of these equivalences is a straightforward manipulation of the definitions, combined with Theorem~\ref{thm:nhbdlat}.

\section{Algorithm for estimating Markov boundaries}
\label{app:gs}

Algorithm~\ref{alg:bdy:comp} below gives the pseudocode for the Grow-Shrink algorithm \citep{margaritis1999bayesian} for estimating Markov boundaries, adapted to the graphoid setting we have adopted.

\begin{algorithm}[H]
            \begin{algorithmic}[1]
                \Function{computeMB}{$j$, $S$, $\graphoid$} 
                \State $m\gets\emptyset$
                \State $w\gets \{i\in S-m : j\notindep i\given m\}$
                \While{ $w \ne \emptyset$ } \Comment{Forward phase}
                \State $m\gets m\cup w$
                \State $w\gets \{i\in S - m: j\notindep i\given m\}$
                \EndWhile
                \ForAll{$i\in m$} \Comment{Backward phase}
                \If {$j\indep i\given m_{-i}$}
                \State $m\gets m_{-i}$
                \EndIf
                \EndFor
                \Return{$m$}
                \EndFunction
            \end{algorithmic}
            \caption{Compute the Markov boundary $\latmin(j;S)$ for a given $j$ and $S \subseteq \ground_{-j}$. \label{alg:bdy:comp}}
\end{algorithm}

\section{Proofs of main results}
\label{sec:proofs:main}

\subsection{Proof of Theorem~\ref{thm:nhbdlat}}
\label{sec:proof:main}

We now prove Theorem~\ref{thm:nhbdlat}. 
First, we establish several simple lemmas that will prove useful. 
As a guide to the reader, Lemmas~\ref{lem:simple:USW} and~\ref{lem:lat:equal} below are intermediate, whereas Lemmas~\ref{lem:reduced:min} and~\ref{lem:inc:equiv} are key lemmas that will be invoked in the main proof.

The following is an immediate consequence of decomposition \ref{defn:gr:decomp}:
\begin{lemma}\label{lem:simple:USW}
    Let $m = \inf\mathfrak{t}_{j}(S)$. Then $j\indep U-m\given m$ for all $U\in\mathfrak{t}_{j}(S)$.
\end{lemma}

\noindent
We will also need the following key technical lemma regarding $\nhbdlatred_{j}(S)$:
\begin{lemma}
\label{lem:reduced:min}
For any $U\in\mathfrak{t}_{j}(S)$, we have $\inf\mathfrak{t}_{j}(U)=\inf\mathfrak{t}_{j}(S)$.
\end{lemma}
\begin{proof}

Let $m_{0}=\inf\mathfrak{t}_{j}(U)$ and $m=\inf\mathfrak{t}_{j}(S)$. By Lemma~\ref{lem:simple:USW}, $j\indep U-m\given m$, whence $m\in\mathfrak{t}_{j}(U)$ by definition. Thus $m_{0}\subseteq m\subseteq U$. Since we have established $m \subseteq U$,
by decomposition \ref{defn:gr:decomp}, $j\indep U-m_{0}\given m_{0}$ implies $j\indep m-m_{0}\given m_{0}$. Combine this with $j \indep S-m \given (m-m_0) \cup m_0$, using contraction~\ref{defn:gr:contr}, to deduce $j\indep S-m_{0}\given m_{0}$. Here, we have used $(S-m) \cup (m-m_0) = S - m_0$ which follows from $m_0 \subseteq m \subseteq S$.
Thus $m_{0}\in\mathfrak{t}_{j}(S)$, and hence $m\subseteq m_{0}$. We conclude that $m=m_{0}$, as desired. 
\end{proof}
\noindent
In the previous lemma, the condition $U\in\mathfrak{t}_{j}(S)$ cannot be removed, as the following example shows.
\begin{example}
In general, if $U \subseteq S$, we will not even have $\inf\mathfrak{t}_{j}(U) \subseteq \inf\mathfrak{t}_{j}(S)$. To see this, consider the set $S_3 = \{4,6,12,14\}$ and $j=3$ from Example~\ref{exa:bdy}, so that $\inf\mathfrak{t}_{j}(S_3) = \{6,12\}$. Take $U = \{4,6,14\} \subseteq S_3$ and note that $U \notin \mathfrak{t}_{j}(S_3) = [\{6,12\}, \{4,6,12,14\}]$. We have $\inf\mathfrak{t}_{j}(U) = \{4,6\}$ which is incomparable to $\inf\mathfrak{t}_{j}(S_3)$.
\end{example}

The next two lemmas, which are useful in their own right, illustrate how elements of $\nhbdlat_{j}(S)$ can be partitioned via knowledge of separation in the underlying graphoid. These results are crucial to property \ref{thm:nhbdlat:sup} in Theorem~\ref{thm:nhbdlat}.

\begin{lemma}
\label{lem:lat:equal}
$T\in\nhbdlat_{j}(S)\implies\nhbdlat_{j}(S)=\nhbdlat_{j}(T)$.
\end{lemma}

\begin{proof}
Since $T\in\nhbdlat_{j}(S)$, we have $\latmin(j;S)=\latmin(j;T)$. Now for any $A\in\nhbdlat_{j}(S)$ we have $\latmin(j;A)=\latmin(j;S)=\latmin(j;T)$, whence $A\in\nhbdlat_{j}(T)$. Thus $\nhbdlat_{j}(S)\subseteq \nhbdlat_{j}(T)$, and similarly $\nhbdlat_{j}(T)\subseteq \nhbdlat_{j}(S)$.
\end{proof}

\begin{lemma}
\label{lem:inc:equiv}
Fix $S\subseteq\ground_{-j}$ and let $A\subseteq\ground_{-j}$ and $B\in\nhbdlat_{j}(S)$ be disjoint. Then $j\indep A\given B$ if and only if $A\cup B\in\nhbdlat_{j}(S)$.
\end{lemma}
\begin{proof}
Since $B\in\nhbdlat_{j}(S)$, we have $\nhbdlat_{j}(B)=\nhbdlat_{j}(S)$ by Lemma~\ref{lem:lat:equal}. Hence $A\cup B\in\nhbdlat_{j}(S)\iff A\cup B\in\nhbdlat_{j}(B)$, which is equivalent to $B\in\nhbdlat_{j}(A\cup B)$. But by definition this means
\begin{align*}
j\indep (A\cup B)-B\given B
&\iff
j\indep A\given B.
\end{align*}
\end{proof}
\noindent
Finally, we proceed with the proof of Theorem~\ref{thm:nhbdlat}.

\begin{proof}[Proof of Theorem~\ref{thm:nhbdlat}]
We break the proof of \ref{thm:nhbdlat:lat} into three parts: (i) Existence of meets, (ii) Existence of joins, and (iii) Convexity. Throughout, we assume that $T,R\in\mathfrak{T}_{j}(S)$. Then, by the definition of $\mathfrak{T}_{j}(S)$, we have $\inf\mathfrak{t}_{j}(T) =  \inf\mathfrak{t}_{j}(R) = m$ for some set $m \subseteq T, R$. It follows that $\mathfrak{t}_{j}(T)=[m,T]$ and $\mathfrak{t}_{j}(R)=[m,R]$; see the discussion following Lemma~\ref{lem:nhbdlatred:lat}.

(i) We wish to show that $T\cap R\in\mathfrak{T}_{j}(S)$. 
Since $m\subseteq T$ and $m\subseteq R$, it follows that $m\subseteq T\cap R\subseteq T$, whence $T\cap R\in\mathfrak{t}_{j}(T)$. Now apply Lemma~\ref{lem:reduced:min} to deduce that $\inf\mathfrak{t}_{j}(T\cap R)=\inf\mathfrak{t}_{j}(T)=m$. 
Thus $T\cap R\in\mathfrak{T}_{j}(S)$.

(ii) We wish to show that $T R\in\mathfrak{T}_{j}(S)$.  Since $T \in \mathfrak{t}_{j}(T)$ and $m = \inf \mathfrak{t}_{j}(T)$, Lemma~\ref{lem:simple:USW} implies $j \indep T-m \given m$. Similarly, $j \indep R-m \given m$. By composition~\ref{defn:gr:comp}, 
$j \indep (TR)-m \given m$,
whence $m \in \mathfrak{t}_{j}(T R)$. Since both $m$ and $TR$ belong to $\mathfrak{t}_{j}(T R)$, and $m \subseteq T \subseteq TR$, it follows from the convexity of $\mathfrak{t}_{j}(T R)$ (Lemma~\ref{lem:nhbdlatred:lat}) that $T \in \mathfrak{t}_{j}(T R)$. Applying Lemma~\ref{lem:reduced:min}, we have $\inf \mathfrak{t}_{j}(T) = \inf \mathfrak{t}_{j}(T R)$, the desired result.

(iii) Suppose $T\subseteq R$ and let $U$ satisfy $T\subseteq U\subseteq R$. Since $m\subseteq T$ we have $U\in[m,R]=\mathfrak{t}_{j}(R)$. Then Lemma~\ref{lem:reduced:min} implies $\inf\mathfrak{t}_{j}(U)=\inf\mathfrak{t}_{j}(R)$, i.e., $U\in\mathfrak{T}_{j}(S)$.

To prove \ref{thm:nhbdlat:bdy}, let $m'=\inf\nhbdlat_{j}(S)$ and note that the definition of the neighbourhood lattice implies that $m'\subseteq\latmin(j;S)$. Now suppose that $m'$ is a proper subset of $\latmin(j;S)$. 
Then since $m'\in\nhbdlat_{j}(S)$, we have $\latmin(j;S)=\latmin(j;m')\subseteq m' \subsetneq \latmin(j;S)$, which is a contradiction. 
The claim \ref{thm:nhbdlat:sup} follows immediately from Lemma~\ref{lem:inc:equiv}.
Finally, to prove \ref{thm:nhbdlat:decomp}, it is enough to observe that $S\sim T\iff\latmin(j;S) = \latmin(j;T)$ defines an equivalence relation on $2^{\ground_{-j}}$.
\end{proof}

\subsection{Proof of Theorem~\ref{thm:poly:N}}
\label{sec:proof:full:decomp}
To reduce notational burden, let us drop the index $j$ and write $\nhbdlatd=\nhbdlatd_j$ through the remainder of this section.
We furthermore write $k:=|\nhbdlatd|$
for the number of lattices in the lattice decomposition $\nhbdlatd$.
The lattice decomposition $\nhbdlatd$ can be computed recursively, as detailed in Algorithm~\ref{alg:lattice:decomp}. The key in this algorithm is step~\ref{step:finding:uncovered:S} where given a partial decomposition, say, $\nhbdlatd'$, one needs to find a set $S$ that is not covered by $\nhbdlatd'$, i.e., $S \notin \bigcup \nhbdlatd'$. If no such set exists, we conclude that $\nhbdlatd'$ covers the power set (i.e., $\bigcup \nhbdlatd' = 2^{\ground_{-j}}$), hence it is in fact the full lattice decomposition and the algorithm terminates. For a general set system $\nhbdlatd'$, performing step~\ref{step:finding:uncovered:S} is NP-hard. However, since in our case the underlying lattices are mutually disjoint, it is possible to produce an uncovered set in polynomial time:

\begin{lemma}[\textproc{findUncoveredSet} Algorithm]\label{lem:decide:set:system}
    Given a set system $\nhbdlatd' =\{\nhbdlat^1,\dots,\nhbdlat^k\}$ where $\nhbdlat^\ell = [\latmin^\ell,\latmax^\ell] \subseteq 2^{\ground_{-j}}$ for all $\ell=1,\dots,k$ and $\nhbdlat^\ell \cap \nhbdlat^{\ell'} = \emptyset$ for $\ell \neq \ell'$, there is a polynomial-time algorithm  to produce a set not covered by $\nhbdlatd'$, 
when such a set exists. 
The algorithm requires at most $O( |\ground|^2 |\nhbdlatd'|)$ set operations, and hence runs in $O(|\ground|^{3}|\nhbdlatd'|)$ time. Moreover, there is a variant of the algorithm that can find all uncovered subsets of size $\le s$ with $O(2^s |\ground| |\nhbdlatd'| )$ set operations.
\end{lemma}

\noindent
Although it is trivial to check whether or not $\nhbdlatd'$ covers $2^{\ground_{-j}}$ (e.g. simply by comparing cardinalities), Lemma~\ref{lem:decide:set:system} goes a step further to produce a certifying set.
There are multiple polynomial-time algorithms that can produce a subset of $\ground_{-j}$ uncovered by set-system $\nhbdlatd'$, as described by Lemma~\ref{lem:decide:set:system}.  We refer to a generic such algorithm as \textproc{findUncoveredSet}$(\nhbdlatd', \ground_{-j})$. A particular instance is given in the proof of the lemma.

\begin{algorithm}[t]
    \begin{algorithmic}[1]
\State Set $\nhbdlat \gets \textproc{computeLattice}(j, \ground_{-j},\graphoid)$.
\State Initialize $\nhbdlatd_j \gets \{\nhbdlat\}$ and $\text{powerSetNotCovered} \gets \text{True}$.
        
        \While{ powerSetNotCovered }
\label{step:finding:uncovered:S}
        \State $S \gets \textproc{findUncoveredSet}(\nhbdlatd_j, \ground_{-j})$.
\If {Step 4 is successful}
        \State $\nhbdlat \gets \textproc{computeLattice}(j, S,\graphoid)$. 
        \State $\nhbdlatd_j \gets \nhbdlatd_j \cup \{\nhbdlat\}$.
\Else
        \State $\text{powerSetNotCovered} \gets \text{False}$.       
        \EndIf
        \EndWhile

    \end{algorithmic}
    \caption{Compute the lattice decomposition $\nhbdlatd_j$ for a given $j$~\label{alg:lattice:decomp}.}
\end{algorithm}

A set operation in the statement of Lemma~\ref{lem:decide:set:system} involves the computation of the size of at most two set differences. We refer to the proof in Section~\ref{sec:proof:comp} for more details.
Since the total number of neighbourhood lattices
found by Algorithm~\ref{alg:lattice:decomp} could potentially be much smaller than $2^{d-1}$ where we recall $d:=|\ground|$, the overall procedure could result in substantial savings relative to the naive approach of checking all the $2^{d-1}$ possible subsets. In general, 
by combining Proposition~\ref{prop:lat:comp:complex} and Lemma~\ref{lem:decide:set:system}, we conclude that Algorithm~\ref{alg:lattice:decomp} achieves the claimed complexity bound in Theorem~\ref{thm:poly:N}.

\subsection{Proof of Theorem~\ref{thm:highd:consist:nongauss}}
\label{app:proof_of_theorem_ref_thm_highd_consist}

Let $a \wedge b := \min\{a,b\}$.
    By a result of~\cite{gcm}, under the conditions of the theorem, for any $i,j \in [d]$ and $S \subseteq [d]$ with $|S| \le t$, we have
    \begin{align}
        \label{eq:rho:concent}
        \pr\bigl( |\gcorrh_n(i,j | S) - \gcorr(i,j|S)| > \gamma \bigr)
        \le C_1 e^{-  2 c_1  n \gamma^2},
    \end{align}
    for any $\gamma \in [\eps_n, c_2 (1 \wedge \sigma^2)]$.
    Let $\mathcal W$ be the collection of all triples $(i,j,S)$ that appear during the run of the population version of Algorithm~\ref{alg:lattice:decomp:sparse} for calculating $\nhbdlatd_j(t)$. We have $|\mathcal W| \lesssim t^3 d^{t+1}$. By union bound, and using $t \le d$,
    \begin{align*}
        \pr\bigl( \max_{(i,j,S) \, \in \, \mathcal W}| \rhoh_n(i,j|S) -\rho(i,j|S)| > \gamma\bigr) \le C_1  d^{t+4} \exp\bigl(-2 c_1 n \gamma^2\bigr).
    \end{align*} 
    Take $\gamma = \alpha/4$. Assuming that $c_1 n \gamma^2 \ge \log(d^{t+4})$, with probability at least $1 - C_1 e^{-c_1 n\gamma^2}$, for all $(i,j,S) \in \mathcal W$, we have
    \begin{align*}
        \begin{cases}
            |\rhoh_n(i,j|S)| < \alpha/4, & \text{if}\; \rho(i,j|S) = 0, \\
            |\rhoh_n(i,j|S)| > 3\alpha/4, & \text{if}\; \rho(i,j|S) \neq  0.
        \end{cases}
    \end{align*} 
    Let us refer to the above event as $\mathcal A$. Then, as long as $\tau \in [\alpha/4, 3\alpha/4]$, on event $\mathcal A$, testing based on $|\rhoh(i,j,S)| \ge \tau$ produces the same result as testing based on $\rho(i,j,S) \neq 0$ for all triples $(i,j,S) \in \mathcal W$. This shows that on event $\mathcal A$, the data-driven Algorithm~\ref{alg:lattice:decomp:sparse} encounters the exact sequence of triples encountered by population-level Algorithm~\ref{alg:lattice:decomp:sparse}, hence producing the exact same output. 
    
Condition $ c_1 n \gamma^2 \ge \log(d^{t+4})$ holds if 
$n \alpha^2 \ge (16/c_1) (t+4)\log d.$
    Similarly, condition~$\gamma \in [\eps_n, c_2 (1 \wedge \sigma^2)]$ holds if 
    $
    4 \eps_n \le \alpha \le 4 c_1 (1 \wedge \sigma^2).
    $
    Putting the pieces together, all the conditions are satisfied if
    \[
\frac{16}{c_1} \cdot (t+4) \frac{ \log d}{n}
 + 16 \eps_n^2  \;\le\; \alpha^2 \;\le\;  4 c_2 (1 \wedge \sigma^2),
    \]
    which gives the desired result.

\subsection{Proof of Theorem~\ref{thm:conn:comp}}
\label{sec:proof:thm:conn:comp}
We now prove Theorem~\ref{thm:conn:comp}.
First, we introduce some new notation and a lemma. Write $\SEP(A,C,B)$ for disjoint sets $A,B$ and $C$ when $C$ separates $A$ and $B$. This means that any path from a node $i \in A$ to node $j \in B$ should pass through a node $k \in C$. This includes the case where there is no path from $A$ to $B$.

\begin{figure}[t]

\centering
    \begin{subfigure}[b]{0.16\textwidth}
        \begin{tikzpicture}
            \begin{scope}[every node/.style={circle,thick,draw, inner sep =.5mm}]
                \node (A) at (0,0) {$A$};
                \node (C_1) at (1,.6) {$C_1$};
                \node (C_2) at (1,-.6) {$C_2$};
                \node (B) at (2,0) {$B$};   
            \end{scope}

            \draw[dashed] (A) -- (C_1);
            \draw[dashed] (A) -- (C_2);
            \draw[thick] (C_1) -- (C_2);
            
        \end{tikzpicture}  
        \caption{} 
    \end{subfigure}
    ~\qquad\qquad~
    \begin{subfigure}[b]{0.16\textwidth}
        \begin{tikzpicture}
            \begin{scope}[every node/.style={circle,thick,draw, inner sep =.5mm}]
                \node (A) at (0,0) {$A$};
                \node (C_1) at (1,.6) {$C_1$};
                \node (C_2) at (1,-.6) {$C_2$};
                \node (B) at (2,0) {$B$};
            \end{scope}
            
            \draw[dashed] (A) -- (C_1);
            \draw[dashed] (A) -- (C_2);
            \draw[dashed] (C_1) -- (B);
        \end{tikzpicture}  
        \caption{} 
        \label{fig:1:a}
    \end{subfigure}

    \caption{Illustration of the reduction property. 
        We break the cases respecting $\SEP(A,C_1 \uplus C_2,B)$ into two groups, depending on whether there is a path between $C_1$ and $C_2$ (a), or  not (b). The dashed lines are possible paths.  In none of the cases there could be a path between $C_2$ and $B$ because of the assumption $\SEP(C_2,A,B)$. The same assumption precludes a path between $C_1$ and $B$ in cases in (a), because of the existence of a path from $C_2$ to $C_1$. As seen from this figure, in all possible cases, $A$ and $B$ are separated by $C_1$.
\label{fig:reduct:property}}\end{figure}

\begin{lemma}\label{lem:add:red}
    Separation has these properties:
    \begin{itemize}
        \item \emph{Additive property}:  $\SEP(A,C_1,B_1)$ and $\SEP(A,C_2,B_2)$ implies $\SEP(A,(C_1\cup C_2),(B_1 \cup B_2))$.
        
\item \emph{Reduction property}: $\SEP(A,C_1 \uplus C_2,B)$ and $\SEP(C_2,A,B)$ implies $\SEP(A,C_1,B)$.
    \end{itemize}
    
\end{lemma}

\begin{proof}
    We only show the reduction property, for which it is enough to show that there is no path from $A$ to $B$ that passes through $C_2$ and not $C_1$. If so, a portion of it gives a path from $C_2$ to $B$ that does not pass through $C_1$, and clearly not $A$ since the paths do not self-cross. But then $A$ does not separate $C_2$ and $B$, a contradiction. Figure~\ref{fig:reduct:property} illustrates the argument.
\end{proof}

\begin{proof}[Proof of Theorem~\ref{thm:conn:comp}]
By induction, we can reduce to the case where the nodes are partitioned into two disjoint 
components  $G_1$ and $G_2$, after the removal of $j$. This implies $\SEP(G_1,j,G_2)$. For any set $A$, let $A_k := A \cap G_k, k =1,2$. We have
    \begin{align}\label{eq:sep:decomp}
        \SEP(j,A,(S\setminus A)) \quad \iff\quad  \SEP(j,A_1,(S_1\setminus A_1)) \quad \text{and}\quad 
        \SEP(j,A_2,(S_2\setminus A_2)).
    \end{align}
    To see this, note that the RHS implies $\SEP(j,A, (S_1\setminus A_1)) \cup (S_2\setminus A_2)$ by the additive property (Lemma~\ref{lem:add:red}). But $(S_1\setminus A_1) \cup (S_2\setminus A_2) = S \setminus A$ since $G_1$ and $G_2$ are disjoint. Now assume the LHS. Then,  clearly $\SEP(j,A,(S_1\setminus A_1))$. But we also have $\SEP(A_2,j,(S_1\setminus A_1))$. Applying the reduction property (Lemma~\ref{lem:add:red}), we get $\SEP(j,A_1,(S_1\setminus A_1))$. The other implication is similar, hence we get the RHS.
    
    Equipped with~\eqref{eq:sep:decomp}, we can prove part~(a). Let $\Ss$ be the smallest subset of $S$ separating $j$ and $S \setminus \Ss$ 
and let $\Ss_1$ and $\Ss_2$ be its components (i.e, $\Ss_k = \Ss \cap G_k,\;\; k = 1,2$). Then, by~\eqref{eq:sep:decomp} $\Ss_k$ separates $j$ and $S_k \setminus\Ss_k$, for $k=1,2$. It remains to show that $\Ss_k$ is the smallest such set (for each $k$). Suppose that there is a proper subset $S'_1$ of $\Ss_k$ that separates $j$ and $S \setminus S'_1$. Then, applying~\eqref{eq:sep:decomp} with $A_1 = S'_1$ and $A_2 = \Ss_2$ and $A = S':= S'_1 \cup \Ss_2$, we conclude that $S'$ separates $j$ and $S\setminus S'$. But $S'$ is a proper subset of $\Ss$, violating the assumption that $\Ss$ is the minimal set with such property. It follows that we should have $m_j(S_1) = \Ss_1$ and $m_j(S_2) = \Ss_2$ which proves part~(a).
    
    For part~(b), let $\Ss_k = m_j(S_k), k =1,2$ and $\Ss = m_j(S)$. From part~(a), we have $\Ss = \Ss_1 \cup \Ss_2$. Fix $r \in \ground_j$. We claim that  
    \begin{align}\label{eq:sep:decomp:2}
        \SEP(r,\Ss,j) \quad \iff \quad \SEP(r,\Ss_1,j) \;\;\text{or}\;\; \SEP(r,\Ss_2,j).
    \end{align}
    The $\Leftarrow$ implication is clear. For the other direction, assume that $r$ is separated from $j$ by $\Ss$, i.e., all the paths from $r$ to $j$ pass through $\Ss$. WLOG, assume that $r \in G_1$ (the case $r \in G_2$ is similar). Suppose that there is a path from $r$ to $j$ that passes through $\Ss_2$. A portion of this path gives a path from $r$ to $\Ss_2$, hence to $G_2$, that does not pass through $j$, contradicting $\SEP(G_1,j,G_2)$. Hence all the paths from $r$ to $j$ should pass through $\Ss_1$, i.e., $\SEP(r,\Ss_1,j)$, proving~\eqref{eq:sep:decomp:2}.
    From~\eqref{eq:sep:decomp:2}, we conclude that $E_j(\Ss) = E_j(\Ss_1) \cup E_j(\Ss_2)$. Combined with characterization of $M_j(S)$ in \ref{grint:latmax},
this proves the first equality in part(b).
    
    To see the second equality in~(b), recall that $\Tc_j(S_k;G_k)$ is the lattice with ground set restricted to $G_k \cup \{j\}$. First, the minimal element of $\Tc_j(S_k;G_k)$ is the same as that of $\Tc_j(S_k)$. This is true since a subset $A$ of $S_k$ separates $j$ and $S_k \setminus A$ in $G_k \cup \{j\}$ iff it does so in $G$. (If the separation happens in $G_k \cup \{j\}$ but not in $G$, there will be a path from $S_k \setminus A$ to $j$ passing through some $G_r$, $r\neq k$ contradicting disjointness of $G_k$ and $G_r$.) Even more directly, a minimal subset of $\Tc_j(S_k)$ is a subset of $S_k$, hence $G_k$, and the restriction in  $\Tc_j(S_k;G_k)$ is automatically satisfied.  Knowing that the minimal element of $\Tc_j(S_k;G_k)$ is $S^*_k$, we invoke \ref{grint:latmax}
restricted to $G_k \cup \{j\}$ to conclude that $M_j(S_k;G_k) = S_k^* \cup E_j(S_k^*;G_k)$ where $E_j(S_k^*;G_k) = \{i \in G_k:\;  \text{$S^*_k$ separates $i$ and $j$}  \}$. Going through the argument leading to~\eqref{eq:sep:decomp:2}, it is clear that we can replace~\eqref{eq:sep:decomp:2} with
    \begin{align}
        \SEP(r,\Ss,j) \quad \iff \quad [\SEP(r,\Ss_1,j) \;\text{and}\;r \in G_1 ] \;\;\text{or}\;\; [\SEP(r,\Ss_2,j)  \;\text{and}\;r \in G_2 ]
    \end{align}
    showing that $E_j(\Ss) = E_j(\Ss_1;G_1) \, \cup\, E_j(\Ss_2;G_2)$. Combined with the expression for $M_j(S_k;G_k)$, the desired result follows.
\end{proof}

\subsection{Proofs for Section~\ref{sec:proj}}
\label{sec:proj:proofs}

This section contains the proofs of Theorems~\ref{thm:reg:lattice} and~\ref{thm:abs:proj}.

\begin{proof}[Proof of Theorem~\ref{thm:reg:lattice}]

    For  $L \in B(\Hil)$,  let $[L]$ be the range projection of $L$, i.e., projection onto the closure of the range of $L$. For any $Q \in \Prj(\Hil)$ and $A  \in  B(\Hil)$, we have $[Q A] \le Q$, since the range of $QA$ is included in the range of $Q$, and the range of $Q$ is closed. We also have the identity $  B = [B] B$ for all $B \in B(\Hil)$.

    Recall 
that $\ker(A^*) = \ran(A)^\perp$ where $A^*$ is the adjoint of $A$. Since $x \in \ran(A)^\perp$ if and only if $[A] x = 0$, we have
    \begin{align}
        \ker(A^*) = \ker([A]), \quad \forall A \in B(\Hil).
    \end{align}
    This in turn is equivalent to 
    \begin{align}\label{eq:adj:range:proj}
        A^* B = 0 \iff  [A] B = 0 \quad \forall B \in B(\Hil).
    \end{align}
    We also need the following result.
    \begin{lemma}\label{lem:additive:proj:join}
        In a finite-dimensional  $\Hil$, we have $P \vee Q = P + [P^\perp Q]$ for all $P,Q \in \Prj(\Hil)$.
    \end{lemma}
    We will use the equivalent characterization~\ref{thm:nhbdlat:equiv:3} from Theorem~\ref{thm:nhbdlat:equiv}. Assume that $U \subseteq S$. Let $P = P_U$, $Q = P_{S-U}$ and $R = P_S$ and note that
    \begin{align}\label{eq:PS:decomp}
        R = P \vee Q = P + [P^\perp Q]
    \end{align}
     where the second equality is by Lemma~\ref{lem:additive:proj:join}. We have that $U \in \Tc_j(S)$ iff $P P_j = R P_j$. Using~\eqref{eq:PS:decomp}, this is equivalent to $[P^\perp Q] P_j = 0$. Applying~\eqref{eq:adj:range:proj}, we equivalently have $(P^{\perp} Q)^* P_j= 0$ which is equivalent to $Q P^{\perp} P_j = 0$, since projections are self-adjoint. Note that this latter statement is just $j \indep S-U \given U$ which is the characterization in~\ref{thm:nhbdlat:equiv:3}, that is, $U \in \nhbdlat_j(S)$.
     
     Now consider $U \supset S$ and redefine $P = P_S$, $Q = P_{U-S}$ and $R = P_U$. The exact same argument as above goes through leading to $j \indep U-S \given S$  which is the characterization in~\ref{thm:nhbdlat:equiv:3}. This completes the proof.
\end{proof}

\begin{proof}[Proof of Lemma~\ref{lem:additive:proj:join}]
Since every $x \in \Hil$ can be decomposed into a component in $\ran(P \vee Q)$ and a component in $\ran(P \vee Q)^\perp$, it is enough to show that the operator equality holds only over these two subspaces.
    
    We first make a couple of observations: We have $\ran(P^\perp Q) \subseteq \ran(P \vee Q)$, that is, $[P^\perp Q] \le P \vee Q$. To see this, let $x \in \ran(P^\perp Q)$. Then, 
    $$x = P^\perp Q z = (I-P) Qz = Qz - PQ z \in \ran(P \vee Q)$$ since it is a linear combination of a component in $\ran(Q)$ and component in $\ran(P)$.
    
    Now assume that $x \in \Hil$ satisfies $(P \vee Q) x = 0$. Multiplying both sides by $P$ and noting that $P \le P \vee Q$, we have $P x = 0$. Similarly, multiplying both sides by $[P^\perp Q]$ and using $[P^\perp Q] \le P \vee Q$, we get $[P^\perp Q]  x = 0$. Thus, $(P + [P^\perp Q])x = 0$ as desired.
    
    Next, assume that  $x \in \Hil$ is such that $(P \vee Q) x = x$. We have $x = x_1 + x_2$ where $x_1 \in \ran(P)$ and $x_2 \in \ran(Q)$. Let $x_3 = x_1 + P x_2 \in \ran(P)$. Then, $x = x_3 + P^\perp x_2$. We can write $x_2 = Q z$ for some $z \in \Hil$. Let $R = [P^\perp Q]$ and note that $R$ and $P$ are orthogonal projections, i.e., $PR = 0$. We have 
    \begin{align*}
        (P + R) x &= (P+R)(x_3 + P^\perp Q z) \\
        &= Px_3 + R P^\perp Q z = x_3 + P^\perp Q z = x,
    \end{align*}
    which is the desired result. The last equality uses $[A]A = A$ which holds for any $A \in B(\Hil)$. The proof is complete.
\end{proof}

\begin{proof}[Proof of Theorem~\ref{thm:abs:proj}]
    To simplify, let $\mathcal Q = \mathcal Q(A, P)$. We will also continue with the notation introduced in the proof of Theorem~\ref{thm:reg:lattice}.
    Now consider a subset $\mathcal R \subseteq \mathcal Q$. For all $Q \in \mathcal R$, we have
    \begin{align}\label{eq:PA:Q}
        [P A] = [Q A] \le Q,
    \end{align}
    hence $[P A]$ is a lower bound on $\mathcal R$. Let $R := \bigwedge \mathcal R$, that is, the infimum of the elements of $\mathcal R$ as projections in $\Prj(\Hil)$. By the definition of infimum, $ [P A] \le R  \le Q$ for all $Q \in \mathcal R$. Then, \begin{align}
        [QA] \le R \le Q
    \end{align}
    for all $Q \in \mathcal R$. By Lemma~\ref{lem:proj:order}, we equivalently have  $R = RQ$ and $[QA] = R [QA]$.  Hence,
    \begin{align}\label{eq:RA:PA}
        RA = R QA  = R[QA] QA  = [QA] QA = QA = PA,
    \end{align}
where the second and fourth equality are by $B = [B] B$. Equation~\eqref{eq:RA:PA} shows that 
    $R \in \mathcal Q$, hence $\mathcal Q$ is closed under infima. 
    
    Since $P A = Q A$ is equivalent to $P^\perp A = Q^\perp A$, it follows that $\mathcal Q^\perp := \{Q^\perp :\; Q \in \mathcal Q\}$ is also closed under infima. But since $Q^\perp \le P^\perp$ iff $P \le Q$, we have that $\mathcal Q$ is closed under suprema, showing that $\mathcal Q$ is a complete sublattice.
    
    For the convexity, assume that $P,Q \in \mathcal Q$ and $U \in \Prj(\Hil)$ such that $P \le U \le Q$. We have 
    $
    PA = U PA = U QA = UA
    $, hence $U \in \mathcal Q$. The proof is complete.
\end{proof}

\section{Additional proofs}
\label{sec:proofs}

\subsection{Proof of Lemma~\ref{lem:nhbdlatred:lat}}
\label{sec:proofs:lem:nhbdlatred:lat}

The result follows from two auxiliary lemmas.

\begin{lemma}
    Assume that the independence model has the weak union property~\ref{defn:gr:wkun}. Then, $\nhbdlatred_{j}(S)$ is closed under unions and is convex. More precisely, if $U \in \nhbdlatred_{j}(S)$ and $U \subseteq T \subseteq S$, then $T \in \nhbdlatred_{j}(S)$.
\end{lemma}
\begin{proof}
    By assumption $j \indep S - U \given U$. Taking $B = S-T$ and $D = T-U$, we have $j \indep BD \given U$. By weak union, $j \indep B \given DU$. But $DU = T$ showing that $T \in \nhbdlatred_{j}(S)$ as desired.
\end{proof}

That $\nhbdlatred_{j}(S)$ is closed under set intersection follows from the following lemma.
\begin{lemma}
    In any graphoid, the following extended intersection property holds:
    \begin{align}\label{defn:gr:ext:int}
        A \indep R' V \mid R W, \quad A \indep R V \mid R' W \implies 
        A \indep R R' V \mid W
    \end{align}
\end{lemma}
We note that the intersection property~\ref{defn:gr:int} corresponds to $V = \emptyset$ in~\eqref{defn:gr:ext:int}.
\begin{proof}
    By decomposition~\ref{defn:gr:decomp}, we can drop $V$ from both assumptions to get $A \indep R' \given R W$ and $A \indep R \given R' W$. By intersection~\ref{defn:gr:int}, $A \indep R R' \given W$. By decomposition~\ref{defn:gr:decomp}, we further have $A \indep R \given W$. Combining this with the assumption $A \indep R' V \given R W$ and applying contraction~\ref{defn:gr:contr} gives the desired result.
\end{proof}
To verify  that $\nhbdlatred_{j}(S)$ is closed under intersections, let $U, U' \in \nhbdlatred_{j}(S)$ and set $W = U \cap U'$, $R = U - W$, $R' = U' - W$ and $V = S - (UU')$. Thus, $\{W,R,R',V\}$ forms a partition of $S$. Note that $U \in \nhbdlatred_{j}(S)$ means $j \indep S -U \given U$ which is equivalent to $j \indep R' V \given R W$. Similarly, $U' \in \nhbdlatred_{j}(S)$ is equivalent to $j \indep R V \given R' W$. The extended intersection property then gives $A \indep R R' V \given W$, that is, $A \indep S - W \given W$. This shows that $W \in \nhbdlatred_{j}(S)$, as desired. This finishes the proof of Lemma~\ref{lem:nhbdlatred:lat}.

\subsection{Proof of Lemma~\ref{lem:lattice:compare}}
\label{sec:proof:comparison}

To prove Lemma~\ref{lem:lattice:compare}, we recall some basic facts given a subset $S\subset\ground_{-j}$:
\begin{enumerate}[label=(A\arabic*)]
\item\label{reg:coef} The regression coefficients are given by $\beta_{j}(S) = \Sigma_{S,S}^{-1}\Sigma_{S,j}^{}$;
\item\label{reg:supp} The projection lattice can be written as $\mathcal{T}_j(S) = \big\{U \subseteq\ground_{-j}: \; \supp(\beta_j(U)) = \supp(\beta_j(S)) \big\}$ (cf. \eqref{eq:lattice:alt});
\item\label{reg:ci} If $X\sim\normalN(0,\Sigma)$ for some $\Sigma\posdef0$ then $X_{j}\indep S-m \given m$ for $m=\supp(\beta_j(S))$. 
\end{enumerate}

\begin{proof}[Proof of Lemma~\ref{lem:lattice:compare}]

(a) This follows by direct calculation using \ref{reg:coef} and \ref{reg:supp}.

(b) Suppose $X\sim\normalN(0,\Sigma)$.
Let $\mathcal{T}_{j}(S)=[m,M]$ and $U\in\mathcal{T}_{j}(S)$. 
Since $m=\supp(\beta_j(S))=\supp(\beta_j(U))$, \ref{reg:ci} implies $X_{j}\indep S-m \given m$, so that $m\in\nhbdlat_{j}(S;\prob)$. Similarly, $m\in\nhbdlat_{j}(U;\prob)$. Since these lattices are disjoint, it follows that $\nhbdlat_{j}(S;\prob)=\nhbdlat_{j}(U;\prob)$. But then $U\in\nhbdlat_{j}(S;\prob)$, as desired.
\end{proof}

\subsection{Proof of Lemma~\ref{lem:decide:set:system}}
\label{sec:proof:comp}

\begin{proof}[Proof of Lemma~\ref{lem:decide:set:system}]
There are multiple algorithms to solve this problem. Here we describe one due to Brendan McKay. 
    To simplify notation let $r := d-1 = |[d]_j|$ and identify $[d]_j$ with $[r]$ by relabeling if need be. 
    Given disjoint lattices $\nhbdlat^\ell = [m^\ell,M^\ell] \in 2^{[r]}, \; \ell=1,\dots,K$, and any set $S \subseteq [r]$, let 
    \begin{align*}
    Q(S) = \sum_{\ell=1}^K Q_\ell(S), \quad Q_\ell(S) = | \{ T \in \nhbdlat^\ell :\; S \subseteq T\} |.
    \end{align*}
    Computing $Q_\ell(S)$ is easy and we consider it an atomic set operation. In particular, if $S \not \subseteq M^\ell$ then $Q_\ell(S) = 0$; otherwise, $Q_\ell(S) = 2^{|(M^\ell \setminus m^\ell) \setminus S|}$. Then computing $Q(S)$ can be done by at most $O(K)$ set operations.
    Note that for any set $S \subseteq [r]$ we have $Q(S) \le 2^{r-|S|}$ with equality if and only if all possible subsets of $[r]$ containing $S$ are present in $\nhbdlatd := \cup_{\ell=1}^K \Tc^\ell$.
    
    Now, if $Q(\emptyset) = 2^{r}$, then power set $2^{[r]}$ is covered. Otherwise, there is $i_1 \in [r]$ such that $Q(\{i_1\}) < 2^{r-1}$ which can be found by trying all the $r$ elements of $[r]$. Then, one can find $\{i_1,i_2\}$ such that $Q(\{i_1,i_2\}) < 2^{r-2}$ by trying all $r$ possible cases for $i_2$. Continuing in this manner, we have either of the following two: 
    \begin{enumerate}
        \item After at most $r-1$ rounds, we have $S \subseteq [d]_j$ such that $Q(S) = 0$. This $S$ has the desired property and the algorithm terminates.
        \item We have $S = \{i_1,\dots,i_{r-1}\}$, with $Q(S) < 2$. Since we are not in Case~1, we should have $Q(S) = 1$. This means that either (a) $S \in \nhbdlatd$ in which case the ground set $[r]$ cannot belong to $\nhbdlatd$,      hence we output $[r]$, or  (b) the ground set $[r] \in \nhbdlatd$ in which case the original $S$ cannot belong to $\nhbdlatd$ and we output $S$.
    \end{enumerate}
    In either case, the overall complexity is at most $O(r^2 K)$ set operations.
    
    For the ``moreover'' part of the claim, 
    first let us call a set $S$ eligible if $Q(S) < 2^{r - |S|}$. We can keep track of all eligible sets in each round of the algorithm.
    Let $\mathcal E_t$ be the collection of all eligible sets $S$ of size $t$, which is available at round $t$ of the algorithm.
    Then, we can compute  $\mathcal E_{t+1}$ by considering $S \cup \{i_{t+1}\}$ for all $S \in \mathcal E_t$ and all possible $i_{t+1}$, and verifying eligibility. This requires $O(|\mathcal E_t| r K)$ operations. Hence, running the algorithm for $t=0,1,\dots,s$, requires $O\big(\sum_{t=0}^s 2^t r K\big) = O(2^s r K)$ as claimed.
\end{proof}

\subsection{Proof of Proposition~\ref{prop:geom:lat:isom}}
\label{sec:proof:projlattice}

\begin{proof}[Proof of Proposition~\ref{prop:geom:lat:isom}]
This proposition follows from a more general result: In a geometric lattice, the sublattice generated by an independent set $\Xc$ is isomorphic to the lattice of all finite subsets of $\Xc$ (see Exercise~7 in Chapter 11 of~\cite{nation1998notes}). We give a more direct proof here. Throughout, $S$ and $T$ will be finite subsets of $\Xc$. First, we note that $S \mapsto P_S$ is one-to-one: Assume that $P_S = P_T$, and let $x \in S$. Then, $x = P_S x = P_T x$, hence, $x$ can be written as a linear combination of the elements of $T$. By linear independence, this is only possible if $x \in T$. Similarly, we conclude that $T \subseteq S$, hence $T = S$.
    
    Next, fix some finite $S,T \subseteq \Xc$. Let $\mathcal S := \Span(S)$ and $\mathcal T := \Span(T)$. We recall that $P_S \wedge P_T$ projects onto $\mathcal S \cap \mathcal T$. We claim that $\mathcal S\cap \mathcal T = \Span(S \cap T)$. Take $x \in \mathcal S\cap \mathcal T$. By writing $x$ once as a linear combination of the elements of $S$ and once as that of the elements of $T$, equating the two, and using linear independence, one concludes that $x$ has a representation as a linear combination of only the elements of $S \cap T$. Hence, $\mathcal S\cap \mathcal T \subseteq \Span(S \cap T)$. The other direction is trivial. Hence, $P_S \wedge P_T = P_{S \cap T}$.
    
    Finally, we need to show that $\Span(\mathcal S \cup \mathcal T) = \Span(S \cup T)$. This is always true (even without independence) and can be verified by considering linear combinations. This shows that $P_S \vee P_T = P_{S \cup T}$, finishing the proof.
\end{proof}

\subsection{Miscellaneous results}\label{sec:details:nbhd:lat}

\begin{lemma}
\textup{(G8)} holds for $\Prj(\Hil)$.
\end{lemma}

\begin{proof}
Let $ P_i = xx^T, P_j = yy^T$ and $P_k = zz^T$ where $x,y$ and $z$ are linearly independent unit norm vectors. Here $x y^T $ denotes the tensor product $ x \otimes y$ and for simplicity write $\ip{x,y} = x^T y$.
By the non-degeneracy assumption, $P_C^\perp z \neq 0$. Let $\tilde z = P_C^\perp z / \norm{P_C^\perp z}$.
Note that $P_C \vee P_k = P_C + \tilde z \tilde z^T$. 
Hence $(P_C \vee P_k)^\perp = P_C^\perp - \tilde z \tilde z^T$. The assumption $P_i (P_C \vee P_k)^\perp P_j = 0$ then implies $P_i P_C^\perp P_j =  P_i \tilde z \tilde z^T P_j$. It follows that
\[
0  = P_i \tilde z \tilde z^T P_j = xy^T (x^T \tilde z \tilde z^T y) 
\]
implying that $x^T \tilde z = 0$ or $\tilde z^T y = 0$. The former assertion is equivalent to $x^T P_C^\perp z = 0$ which in turn is equivalent to $P_i P_C^\perp P_k = 0$. The second assertion is similarly equivalent to $P_j P_C^\perp P_k = 0$, finishing the proof.
\end{proof}

\begin{lemma}\label{lem:equiv:lattice:defs}
    Representations~\eqref{eq:lattice:proj:def} and~\eqref{eq:lattice:alt} are equivalent.
\end{lemma}

\begin{proof}Let $\Tc_j(S)$ be as defined in~\eqref{eq:lattice:proj:def}.
    To see the first equality in~\eqref{eq:lattice:alt}, note that $\beta_j(S)^\tpose \dvec = P_S \dvec_j$. Since the covariance matrix of $\dvec$ is positive definite, $\beta_j(S) = \beta_j(T)$ is equivalent to $\beta_j(S)^\tpose \dvec = \beta_j(T)^\tpose \dvec$, and hence to $P_S \dvec_j = P_T \dvec_j$. The result follows since for any operator $L$, $L \dvec_j = 0$ is equivalent to $L P_j = 0$. (We have $L P_j y = \ip{\gv_j,y} L \gv_j$ for all $y$, assuming $\norm{\gv_j}_{\Hil}=1$.) The second equality in~\eqref{eq:lattice:alt} follows from the first equality and positive definiteness of the covariance matrix. 
\end{proof}

\vskip 0.2in
\bibliography{lattice}
\bibliographystyle{abbrvnat}

\end{document}